\documentclass[12pt,a4paper,draft]{amsart} 

\usepackage{amssymb}
\usepackage[initials]{amsrefs}
\usepackage[all]{xy}
\usepackage{ifthen}

\newcommand{\bd}{\mathbf{d}}
\newcommand{\be}{\mathbf{e}}
\newcommand{\bh}{\mathbf{h}}
\newcommand{\bm}{\mathbf{m}}
\newcommand{\bP}{\mathbf{P}}
\newcommand{\bQ}{\mathbf{Q}}
\newcommand{\bR}{\mathbf{R}}
\newcommand{\bX}{\mathbf{X}}
\newcommand{\bY}{\mathbf{Y}}
\newcommand{\blambda}{\mbox{\boldmath $\lambda$}}

\newcommand{\calE}{\mathcal{E}}
\newcommand{\calP}{\mathcal{P}}
\newcommand{\calQ}{\mathcal{Q}}
\newcommand{\calR}{\mathcal{R}}
\newcommand{\calX}{\mathcal{X}}
\newcommand{\calY}{\mathcal{Y}}
\newcommand{\calU}{\mathcal{U}}
\newcommand{\calV}{\mathcal{V}}
\newcommand{\calZ}{\mathcal{Z}}

\newcommand{\bbA}{\mathbb{A}}
\newcommand{\bbM}{\mathbb{M}}
\newcommand{\bbN}{\mathbb{N}}
\newcommand{\bbP}{\mathbb{P}}
\newcommand{\bbV}{\mathbb{V}}
\newcommand{\bbZ}{\mathbb{Z}}

\newcommand{\frakp}{\mathfrak{p}}
\newcommand{\frakO}{\mathfrak{O}}

\def\ffrac#1#2{{\textstyle\frac{#1}{#2}}}

\newcommand{\ts}[1]{{\textstyle#1}}

\let\Im=\undefined
\let\mod=\undefined
\DeclareMathOperator{\GL}{GL} %
\DeclareMathOperator{\id}{id} %
\DeclareMathOperator{\Id}{Id} %
\DeclareMathOperator{\Im}{Im} %
\DeclareMathOperator{\pd}{pd} %
\DeclareMathOperator{\ext}{ext} %
\DeclareMathOperator{\Ext}{Ext} %
\DeclareMathOperator{\Hom}{Hom} %
\DeclareMathOperator{\Ker}{Ker} %
\DeclareMathOperator{\mod}{mod} %
\DeclareMathOperator{\gldim}{gl.dim} %
\DeclareMathOperator{\bdim}{\mathbf{dim}} %

\newtheorem*{coro}{Corollary}
\newtheorem*{lemm}{Lemma}
\newtheorem*{prop}{Proposition}

\newtheorem{theon}{Theorem}[section]

\numberwithin{equation}{subsection}

\title{Geometry of regular modules over canonical algebras}

\author{Grzegorz Bobi\'nski}

\address{Faculty of Mathematics and Computer Science \\ Nicolaus
Copernicus University \\ ul.~Chopina 12/18 \\ 87-100 Toru\'n \\
Poland}

\email{gregbob@mat.uni.torun.pl}


\dedicatory{Dedicated to the memory of Professor Stanis\l aw
Balcerzyk}


\begin{document}

\begin{abstract}
We classify canonical algebras such that for every dimension
vector of a regular module the corresponding module variety is
normal (respectively, a complete intersection). We also prove that
for the dimension vectors of regular modules normality is
equivalent to irreducibility.
\end{abstract}

\maketitle

\makeatletter
\def\@secnumfont{\mdseries} 
\makeatother

\section{Introduction and main result} \label{sectintro}

Throughout the paper $k$ is a fixed algebraically closed field. By
an algebra we always mean a finite dimensional algebra over $k$
and by a module a finite dimensional left module.

In~\cite{Ri}*{3.7} Ringel introduced a class of so-called
canonical algebras (see~\ref{can_def} for a definition). A
canonical algebra $\Lambda$ depends on a sequence $(m_1, \ldots,
m_n)$, $n > 2$, of positive integers greater than $1$, and on a
sequence $(\lambda_3, \ldots, \lambda_n)$ of pairwise distinct
nonzero elements of $k$. In the above situation we say that
$\Lambda$ is a canonical algebra of type $(m_1, \ldots, m_n)$.
These algebras play a prominent role in the representation theory
of algebras. For example their module categories serve as model
categories for module categories of algebras admitting separating
tubular families (see~\cites{LePe, Sk}). The module categories of
canonical algebras are derived equivalent to the categories of
coherent sheaves over weighted projective lines (see~\cite{GeLe}).
Moreover, according to~\cite{Ha}*{Theorem~3.1} every quasi-titled
algebra is derived equivalent either to a hereditary algebra or to
a canonical one.

An important and interesting direction of research in the
representation theory of algebras is study of varieties
$\mod_\Lambda (\bd)$ of $\Lambda$-modules of dimension vector
$\bd$ (see~\ref{var_def}), where $\bd$ is an element of the
Grothendieck group $K_0 (\Lambda)$ (for some reviews of results
see for example~\cites{Bon3, Ge, Kr2}). In particular, varieties
of modules over canonical algebras have been studied.
In~\cite{BobSk1} Skowro\'nski and the author proved that if
$\Lambda$ is a tame canonical algebra and $\bd$ is the dimension
vector of an indecomposable $\Lambda$-module, then $\mod_\Lambda
(\bd)$ is a complete intersection with at most~$2$ irreducible
components. It was also shown that in the above case
irreducibility of $\mod_\Lambda (\bd)$ is equivalent to normality.

For a canonical algebra $\Lambda$ one may distinguish so-called
regular modules (see~\ref{reg_def}). This class of modules also
received special attention from a geometric point of view.
Skowro\'nski and the author showed in~\cite{BobSk2} that if $\bd$
is the dimension vector of a regular module over a tame canonical
algebra $\Lambda$, then the corresponding variety is an
irreducible and normal complete intersection. Similar results for
special cases of wild canonical algebras were obtained by Barot
and Schr\"oer in~\cite{BaSc}. It is also worth mentioning that if
$\bd$ is the dimension vector of a regular module over a canonical
algebra, then descriptions of the semi-invariants with respect to
the natural action of $\GL (\bd)$ were given independently by
Skowro\'nski and Weyman in~\cite{SkWe} and Domokos and Lenzing
in~\cites{DoLe1, DoLe2}.

Our first theorem generalizes to regular modules over arbitrary
ca\-nonical algebra a result obtained for indecomposable modules
over tame canonical algebra in~\cite{BobSk1}.

\begin{theon} \label{theoirr}
Let $\Lambda$ be a canonical algebra and let $\bd$ be the
dimension vector of a regular $\Lambda$-module. Then $\mod_\Lambda
(\bd)$ is normal if and only if it is irreducible.
\end{theon}

Let $a (\bd) = \dim \GL (\bd) - \langle \bd, \bd \rangle$ for $\bd
\in K_0 (\Lambda)$, where $\GL (\bd)$ is the corresponding product
of general linear groups (see~\ref{group_def}) and $\langle -, -
\rangle : K_0 (\Lambda) \times K_0 (\Lambda) \to \bbZ$ is the
Ringel bilinear form (see~\ref{form_def}). We have the following
criterion for a complete intersection.

\begin{theon} \label{theocomp}
Let $\Lambda$ be a canonical algebra and let $\bd$ be the
dimension vector of a regular $\Lambda$-module. Then $\mod_\Lambda
(\bd)$ is a complete intersection if and only if $\dim
\mod_\Lambda (\bd) = a (\bd)$.
\end{theon}

In Propositions~\ref{criterioncomplete} and~\ref{criterion} we
show how the above theorems can be translated into numeric
properties of the Ringel form.

Our aim in this paper is to classify canonical algebras such that
the corresponding module varieties have ``good'' geometric
properties for all dimension vectors of regular modules. It is
done in the following theorem.

\begin{theon} \label{theo1}
Let $\Lambda$ be a canonical algebra of type $(m_1, \ldots, m_n)$.
\begin{enumerate}

\item
The varieties $\mod_\Lambda (\bd)$ are complete intersections for
all dimension vectors $\bd$ of regular $\Lambda$-modules if and
only if
\[
\ffrac{1}{m_1 - 1} + \cdots + \ffrac{1}{m_n - 1} \geq 2 n - 5.
\]

\item
The varieties $\mod_\Lambda (\bd)$ are normal for all dimension
vectors $\bd$ of regular $\Lambda$-modules if and only if
\[
\ffrac{1}{m_1 - 1} + \cdots + \ffrac{1}{m_n - 1} > 2 n - 5.
\]

\end{enumerate}
\end{theon}

Recall that if $\Lambda$ is a canonical algebra of type $(m_1,
\ldots, m_n)$, then $\Lambda$ is of tame (respectively, domestic)
representation type if and only if
\[
\ffrac{1}{m_1} + \cdots + \ffrac{1}{m_n} \geq n - 2 \; (> n - 2).
\]

A natural assumption when dealing with geometric problems is that
$\bd$ is the dimension vector of a sincere module $M$ (i.e., every
simple module occurs as a composition factor of $M$). Such
dimension vectors are also called sincere. We have the
corresponding result in this case.

\begin{theon} \label{theo2}
Let $\Lambda$ be a canonical algebra of type $(m_1, \ldots, m_n)$.
\begin{enumerate}

\item
The varieties $\mod_\Lambda (\bd)$ are complete intersections for
all dimension vectors $\bd$ of sincere regular $\Lambda$-modules
if and only if
\[
\ffrac{1}{m_1 - 1} + \cdots + \ffrac{1}{m_n - 1} \geq 2 n - 5.
\]

\item
The varieties $\mod_\Lambda (\bd)$ are normal for all dimension
vectors $\bd$ of sincere regular $\Lambda$-modules if and only if
either
\[
\ffrac{1}{m_1 - 1} + \cdots + \ffrac{1}{m_n - 1} > 2 n - 5,
\]
or $n = 5$ and $m_i = 2$ for all $i = 1, \ldots, 5$.
\end{enumerate}
\end{theon}

Let $\bh$ be the dimension vector of the multiplicity free sincere
semi-simple $\Lambda$-module (see~\ref{h_def}). It can be observed
from~\cites{BobSk2, DoLe2, SkWe}, that results about $\mod_\Lambda
(\bd)$ depend on whether there exists a regular $\Lambda$-module
$M$ of dimension vector $\bd$ which has a direct summand of
dimension vector $\bh$. Let $\bR'$ be the set of all such
dimension vectors (see also~\ref{PRQ}).

\begin{theon} \label{theo3}
Let $\Lambda$ be a canonical algebra of type $(m_1, \ldots, m_n)$.
\begin{enumerate}

\item
The varieties $\mod_\Lambda (\bd)$ are complete intersections for
all dimension vectors $\bd \in \bR'$ if and only if
\[
\ffrac{1}{m_1 - 1} + \cdots + \ffrac{1}{m_n - 1} \geq 2 n - 5.
\]

\item
The varieties $\mod_\Lambda (\bd)$ are normal for all dimension
vectors $\bd \in \bR'$ if and only if
\[
\ffrac{1}{m_1 - 1} + \cdots + \ffrac{1}{m_n - 1} \geq 2 n - 5.
\]

\end{enumerate}
\end{theon}

We exclude from our considerations the case of canonical algebras
of type $(m_1, m_2)$, since in this case the module varieties are
just affine spaces. However, the above theorems are trivially
satisfied also in this case, if we set $\frac{1}{0} = \infty$.

The paper is organized as follows. In Section~\ref{sectcan} we
present necessary facts about canonical algebras. In
Section~\ref{sectgeom} we collect some useful facts about
varieties of modules, while in Section~\ref{sectproof} we prove
Theorems~\ref{theoirr} and~\ref{theocomp}, and show how to reduce
the proofs of Theorems~\ref{theo1}, \ref{theo2} and~\ref{theo3} to
questions about properties of the Ringel form. Next in
Section~\ref{sectineq} we prove inequalities which show that for
canonical algebras satisfying the conditions of
Theorems~\ref{theo1}, \ref{theo2} and~\ref{theo3}, the
corresponding module varieties have the required properties. On
the other hand, in Section~\ref{sectexamp} we present examples
showing that the above statements do not hold for the remaining
canonical algebras.

The results presented in this paper were obtained while the author
held a one year post-doc position at the University of Bern. The
author gratefully acknowledges the support from the
Schweizerischer Nationalfonds and the Polish Scientific Grant KBN
No.~1 P03A 018 27. The author also expresses his gratitude to
Professor Riedtmann for discussions, which were an inspiration for
this research.

\makeatletter
\def\@secnumfont{\mdseries} 
\makeatother

\section{Facts about canonical algebras} \label{sectcan}

Throughout the paper, by $\bbN$ and $\bbZ$ we denote the sets of
nonnegative integers and integers, respectively. If $i, j \in
\bbZ$, then $[i, j]$ denotes the set of all $l \in \bbZ$ such that
$i \leq l \leq j$.

\subsection{}
Recall that by a quiver $\Delta$ we mean a finite set $\Delta_0$
of vertices and a finite set $\Delta_1$ of arrows together with
two maps $s, t : \Delta_1 \to \Delta_0$, which assign to an arrow
$\gamma \in \Delta_1$ its starting and terminating vertex,
respectively. By a path of length $m \geq 1$ in $\Delta$ we mean a
sequence $\sigma = \gamma_1 \cdots \gamma_m$ of arrows such that
$s \gamma_i = t \gamma_{i + 1}$ for $i \in [1, m - 1]$. We write
$s \sigma$ and $t \sigma$ for $s \gamma_m$ and $t \gamma_1$,
respectively. For each vertex $x$ of $\Delta$ we introduce a path
$x$ of length $0$ such that $s x = x = t x$. We only consider
quivers without oriented cycles, i.e., we assume that there exists
no path $\sigma$ of positive length such that $t \sigma = s
\sigma$.

With a quiver $\Delta$ we associate its path algebra $k \Delta$,
which as a $k$-vector space has a basis formed by all paths in
$\Delta$ and whose multiplication is induced by the composition of
paths. By a relation $\rho$ in $\Delta$ we mean a linear
combination of paths of length at least $2$ with the same starting
and terminating vertex. This common starting vertex is denoted by
$s \rho$ and the common terminating vertex by $t \rho$. A set $R$
of relations is called minimal if for every $\rho \in R$, $\rho$
does not belong to the ideal $\langle R \setminus \{ \rho \}
\rangle$ of $k \Delta$ generated by $R \setminus \{ \rho \}$. A
pair $(\Delta, R)$ consisting of a quiver $\Delta$ and a minimal
set of relations $R$ is called a bound quiver. If $(\Delta, R)$ is
a bound quiver, then the algebra $k \Delta / \langle R \rangle$ is
called the path algebra of $(\Delta, R)$.

\subsection{}
Let $\Lambda$ be the path algebra of a bound quiver $(\Delta, R)$.
It is known that the category $\mod_\Lambda$ of $\Lambda$-modules
is equivalent to the category of representations of $(\Delta, R)$
(see for example~\cite{Ri}*{2.1}). Recall, that by a
representation of $(\Delta, R)$ we mean a collection $(M_x,
M_\gamma)_{x \in \Delta_0, \, \gamma \in \Delta_1}$ of finite
dimensional $k$-vector spaces $M_x$, $x \in \Delta_0$, and
$k$-linear maps $M_\gamma : M_{s \gamma} \to M_{t \gamma}$,
$\gamma \in \Delta_1$, such that $M_\rho = 0$ for all $\rho \in
R$. Here, if $\sigma = \gamma_1 \cdots \gamma_m$ is a path in
$\Delta$, then we write $M_\sigma = M_{\gamma_1} \cdots
M_{\gamma_m}$, and if $\rho = \lambda_1 \sigma_1 + \cdots +
\lambda_n \sigma_n$ is a relation in $\Delta$, then $M_\rho =
\lambda_1 M_{\sigma_1} + \cdots + \lambda_n M_{\sigma_n}$. If $M$
and $N$ are two representations of $(\Delta, R)$, then by a
morphism $f : M \to N$ we mean a collection $(f_x)_{x \in
\Delta_0}$ of linear maps $f_x : M_x \to N_x$, $x \in \Delta_0$,
such that $f_{t \gamma} M_\gamma = N_\gamma f_{s \gamma}$ for all
$\gamma \in \Delta_1$. From now on we identify $\Lambda$-modules
with representations of $(\Delta, R)$. In particular, for each
$\Lambda$-module $M$ we define its dimension vector $\bdim M \in
\bbN^{\Delta_0}$ by $(\bdim M)_x = \dim_k M_x$, $x \in \Delta_0$.

\subsection{} \label{form_def} %
Let $\Lambda$ be the path algebra of a bound quiver $(\Delta, R)$.
For a vertex $x$ of $\Delta_0$ we denote by $\be_x$ the element of
the canonical basis of $\bbZ^{\Delta_0}$ corresponding to $x$. For
$\bd \in \bbZ^{\Delta_0}$ we write $\bd = \sum_{x \in \Delta_0}
d_x \be_x$. Assume that $\gldim \Lambda \leq 2$. We have the
Ringel bilinear form $\langle -, - \rangle : \bbZ^{\Delta_0}
\times \bbZ^{\Delta_0} \to \bbZ$ defined by
\[
\langle \bd', \bd'' \rangle = \sum_{x \in \Delta_0} d_x' d_x'' -
\sum_{\gamma \in \Delta_1} d_{s \gamma}' d_{t \gamma}'' +
\sum_{\rho \in R} d_{s \rho}' d_{t \rho}''.
\]
It is known (see~\cite{Bon1}*{2.2}), that if $M$ and $N$ are
$\Lambda$-modules, then
\[
\langle \bdim M, \bdim N \rangle = [M, N] - [M, N]^1 + [M, N]^2,
\]
where following Bongartz~\cite{Bon2} we write $[M, N] = \dim_k
\Hom_\Lambda (M, N)$, $[M, N]^1 = \dim_k \Ext_\Lambda^1 (M, N)$
and $[M, N]^2 = \dim_k \Ext_\Lambda^2 (M, N)$.

\subsection{} \label{can_def}
Let $\bm = (m_1, \ldots, m_n)$, $n \geq 3$, be a sequence of
integers greater than $1$ and let $\blambda = (\lambda_3, \ldots,
\lambda_n)$ be a sequence of pairwise distinct nonzero elements of
$k$. We define $\Lambda (\bm, \blambda)$ as the path algebra of
the quiver $\Delta (\bm)$
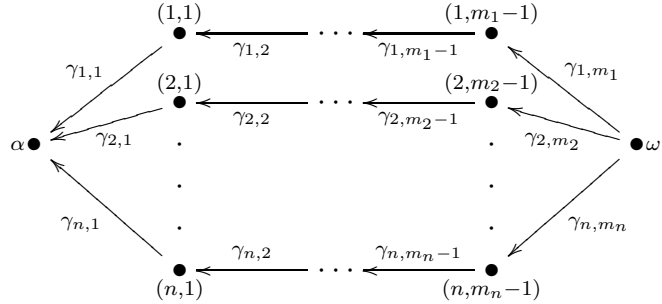
\[
\xymatrix@R=0.25\baselineskip@C=3\baselineskip{%
& \bullet \save*+!D{\scriptstyle (1, 1)} \restore
\ar[lddd]_{\gamma_{1, 1}} & \cdots \ar[l]^-{\gamma_{1, 2}} &
\bullet \save*+!D{\scriptstyle (1, m_1 - 1)} \restore
\ar[l]^-{\gamma_{1, m_1 - 1}}
\\ \\ %
& \bullet \save*+!D{\scriptstyle (2, 1)} \restore
\ar[ld]^{\gamma_{2, 1}} & \cdots \ar[l]^-{\gamma_{2, 2}} & \bullet
\save*+!D{\scriptstyle (2, m_2 - 1)} \restore \ar[l]^-{\gamma_{2,
m_2 - 1}}
\\ %
\bullet \save*+!R{\scriptstyle \alpha} \restore & \cdot & & \cdot
& \bullet \save*+!L{\scriptstyle \omega} \restore
\ar[luuu]_{\gamma_{1, m_1}} \ar[lu]^{\gamma_{2, m_2}}
\ar[lddd]^{\gamma_{n, m_n}}
\\ %
& \cdot & & \cdot
\\ %
& \cdot & & \cdot
\\ %
& \bullet \save*+!U{\scriptstyle (n, 1)} \restore
\ar[luuu]^{\gamma_{n, 1}} & \cdots \ar[l]_-{\gamma_{n, 2}} &
\bullet \save*+!U{\scriptstyle (n, m_n - 1)} \restore
\ar[l]_-{\gamma_{n, m_n - 1}} }
\]
bound by relations
\[
\gamma_{1, 1} \cdots \gamma_{1, m_1} + \lambda_i \gamma_{2, 1}
\cdots \gamma_{2, m_2} - \gamma_{i, 1} \cdots \gamma_{i, m_i}, \,
i \in [3, n].
\]
The algebras of the above form are called canonical. In
particular, we call $\Lambda (\bm, \blambda)$ a canonical algebra
of type $\bm$. It is well known (see for
example~\cite{HaReSm}*{III.4}) that $\gldim \Lambda (\bm,
\blambda) = 2$. If $\bm$ and $\blambda$ are fixed, then we usually
write $\Lambda$ and $\Delta$ instead of $\Lambda (\bm, \blambda)$
and $\Delta (\bm)$, respectively. From now till the end of the
section we assume that $\Lambda = \Lambda (\bm, \blambda)$ is a
fixed canonical algebra.

\label{h_def} %
We write $\be_{i, j}$ instead of $\be_{(i, j)}$ for $i \in [1, n]$
and $j \in [1, m_i - 1]$. For future convenience for $i \in [1,
n]$ by $(i, 0)$ and $(i, m_i)$ we mean $\alpha$ and $\omega$,
respectively. Moreover, if $i \in [1, n]$ and $j \in [0, m_i]$,
then we write $d_{i, j}$ instead of $d_{(i, j)}$. Let $\bh =
\sum_{x \in \Delta_0} \be_x$. For $i \in [1, n]$, we put $\be_{i,
0} = \be_{i, m_i} = \bh - \sum_{j \in [1, m_i - 1]} \be_{i, j}$.
Note that
\begin{align*}
\langle \be_{i, j}, \bd \rangle & = d_{i, j} - d_{i, j - 1}, \, i
\in [1, n], \, j \in [1, m_i],
\\ %
\intertext{and} %
\langle \bd, \be_{i, j} \rangle & = d_{i, j} - d_{i, j + 1}, \, i
\in [1, n], \, j \in [0, m_i - 1],
\end{align*}
and consequently $\langle \bh, \bd \rangle = d_\omega - d_\alpha =
- \langle \bd, \bh \rangle$, for all $\bd \in \bbZ^{\Delta_0}$.

\subsection{} \label{reg_def}
Let $\calP$ ($\calR$, $\calQ$, respectively) be the subcategory of
all $\Lambda$-modules which are direct sums of indecomposable
$\Lambda$-modules $X$ such that
\[
\langle \bdim X, \bh \rangle > 0 \quad (\langle \bdim X, \bh
\rangle = 0, \, \langle \bdim X, \bh \rangle < 0, \text{
respectively}).
\]
The $\Lambda$-modules belonging to $\calR$ are called regular. We
have the following properties of the above decomposition of
$\mod_\Lambda$ (see~\cite{Ri}*{3.7}).

First, $[N, M] = 0$ and $[M, N]^1 = 0$, if either $N \in \calR
\vee \calQ$ and $M \in \calP$, or $N \in \calQ$ and $M \in \calP
\vee \calR$. Here, for two subcategories $\calX$ and $\calY$ of
$\mod_\Lambda$, we denote by $\calX \vee \calY$ the additive
closure of their union. Secondly, $\calR$ decomposes into a
$\bbP^1 (k)$-family $\coprod_{\lambda \in \bbP^1 (k)}
\calR_\lambda$ of uniserial categories. If $\lambda \in \bbP^1 (k)
\setminus \{ \lambda_1, \ldots, \lambda_n \}$, where $\lambda_1 =
0$ and $\lambda_2 = \infty$, then there is a unique simple object
in $\calR_\lambda$ and its dimension vector is $\bh$. On the other
hand, if $\lambda = \lambda_i$ for $i \in [1, n]$, then there are
$m_i$ simple objects in $\calR_{\lambda_i}$ and their dimension
vectors are $\be_{i, j}$, $j \in [1, m_i]$. Finally, one knows
that $\pd_\Lambda M \leq 1$ for $M \in \calP \vee \calR$ and
$\id_\Lambda N \leq 1$ for $N \in \calR \vee \calQ$.

\subsection{} \label{PRQ} %
We denote by $\bP$, $\bR$ and $\bQ$ the sets of the dimension
vectors of the $\Lambda$-modules belonging to $\calP$, $\calR$ and
$\calQ$, respectively. Note that $\bd \in \bR$ if and only if
\[
\bd = p \bh + \sum_{i \in [1, n]} \sum_{j \in [1, m_i]} p_{i, j}
\be_{i, j}
\]
for some nonnegative integers $p$ and $p_{i, j}$, $i \in [1, n]$,
$j \in [1, m_i]$. We know from~\cite{Ri}*{3.7}, that if $\bd \in
\bP$, $\bd \neq 0$, then $d_\alpha > d_\omega \geq 0$ and $d_{i,
j} \geq d_{i, j + 1}$ for all $i \in [1, n]$ and $j \in [0, m_i -
1]$. We show now the converse. Let $\bd$ be as above. Fix
$\lambda_0 \in \bbP^1 (k) \setminus \{ \lambda_1, \ldots,
\lambda_n \}$. It is easy to see that there exists $M \in \calP
\vee \calR_{\lambda_0}$ of dimension vector $\bd$. Indeed, it is
enough to write $\bd = \bd' + \bd''$, where $\bd'' = d_\omega
\bh$. Then obviously there is $M'' \in \calR_{\lambda_0}$ of
dimension vector $\bd''$, and one can easily construct $M' \in
\calP$ of dimension vector $\bd'$, since $d_\omega' = 0$. Since
$[N', N''] = 0$ for $N' \in \coprod_{\lambda \neq \lambda_0}
\calR_\lambda \vee \calQ$ and $N'' \in \calP \vee
\calR_{\lambda_0}$, it follows that $\calP \vee \calR_{\lambda_0}$
is extension closed. In particular, if we assume that the
dimension of the endomorphism ring of $M$ is minimal possible,
then $M = M' \oplus M''$ for $M' \in \calP$ and $M'' \in
\calR_{\lambda_0}$ such that $[M'', M']^1 = 0$ (see for
example~\cite{Ri}*{2.3}). On the other hand, $M'' = p \bh$ for a
nonnegative integer $p$, and $[M'', M']^1 = -\langle p \bh, \bd -
p \bh \rangle = p (d_\alpha - d_\omega) > 0$, if $p > 0$. Thus $p
= 0$, $M'' = 0$, and $M = M' \in \calP$.

Dually, $\bd \in \bQ$, $\bd \neq 0$, if and only if $0 \leq
d_\alpha < d_\omega$ and $d_{i, j - 1} \leq d_{i, j}$ for all $i
\in [1, n]$ and $j \in [1, m_i]$. Thus, each $\bd \in \bQ$ can be
written in a form
\[
\bd = p \bh + \sum_{i \in [1, n]} \sum_{j \in [1, m_i - 1]} p_{i,
j} \be_{i, j} + p_\omega \be_\omega
\]
for some nonnegative integers $p$, $p_\omega$ and $p_{i, j}$, $i
\in [1, n]$, $j \in [1, m_i - 1]$. Consequently, $\bd \in \bR +
\bQ$ if and only if
\[
\bd = p \bh + \sum_{i \in [1, n]} \sum_{j \in [1, m_i]} p_{i, j}
\be_{i, j} + p_\omega \be_\omega
\]
for some nonnegative integers $p$, $p_\omega$ and $p_{i, j}$, $i
\in [1, n]$, $j \in [1, m_i]$. In particular, if $\bd \in \bR +
\bQ$, then there exists a unique presentation
\[
\bd = p^{\bd} \bh + \sum_{i \in [1, n]} \sum_{j \in [1, m_i]}
p_{i, j}^{\bd} \be_{i, j} + p_\omega^{\bd} \be_\omega
\]
such that $p^{\bd}$, $p_\omega^{\bd}$ and $p_{i, j}^{\bd}$, $i \in
[1, n]$, $j \in [1, m_i]$ are nonnegative integers, and for each
$i \in [1, n]$ there exists $j \in [1, m_i]$ such that $p_{i,
j}^{\bd} = 0$. Note that $\bd \in \bR$ if and only if
$p_\omega^{\bd} = 0$. Moreover, $\bd \in \bR'$ if and only if
$p_\omega ^{\bd} = 0$ and $p^{\bd} \neq 0$. Recall that by $\bR'$
we denote the set of all dimension vectors of regular
$\Lambda$-modules which have a direct summand of dimension vector
$\bh$.

\makeatletter
\def\@secnumfont{\mdseries} 
\makeatother

\section{Varieties of modules} \label{sectgeom}

Throughout this section $\Lambda$ is the path algebra of a bound
quiver $(\Delta, R)$ of global dimension at most $2$.

\subsection{}  \label{var_def} %
For $\bd', \bd'' \in \bbN^{\Delta_0}$, let $\bbA (\bd', \bd'') =
\prod_{\gamma \in \Delta_1} \bbM (d_{t \gamma}', d_{s \gamma}'')$,
where by $\bbM (p, q)$ we denote the space of $p \times
q$-matrices with coefficients in $k$. For a dimension vector $\bd
\in \bbN^{\Delta_0}$, $M \in \bbA (\bd, \bd)$ and a path $\sigma =
\gamma_1 \cdots \gamma_m$ of positive length, we put $M_\sigma =
M_{\gamma_1} \cdots M_{\gamma_m}$. We extend this notation to
relations in the standard way. We denote by $\mod_\Lambda (\bd)$
the set of all $M \in \bbA (\bd, \bd)$, such that $M_\rho = 0$ for
all $\rho \in R$. Obviously, $\mod_\Lambda (\bd)$ is an affine
variety. Note that every point $M$ of $\mod_\Lambda (\bd)$
determines a $\Lambda$-module of dimension vector $\bd$ (by taking
$M_x = k^{d_x}$ for $x \in \Delta_0$), which we also denote by
$M$, and every $\Lambda$-module of dimension vector $\bd$ is
isomorphic to $M$ for some $M \in \mod_\Lambda (\bd)$. We call
$\mod_\Lambda (\bd)$ the variety of $\Lambda$-modules of dimension
vector $\bd$. Note that $a (\bd)$ defined in
Section~\ref{sectintro} can be calculated as follows
\[
a (\bd) = \sum_{\gamma \in \Delta_1} d_{s \gamma} d_{t \gamma} -
\sum_{\rho \in R} d_{s \rho} d_{t \rho},
\]
hence $a (\bd)$ is just the dimension of $\bbA (\bd, \bd)$ minus
the number of equations defining $\mod_\Lambda (\bd)$. In
particular, the dimension of each irreducible component of
$\mod_\Lambda (\bd)$ is at least $a (\bd)$.

\label{group_def} %
The product $\GL (\bd) = \prod_{x \in \Delta_0} \GL (d_x)$ of
general linear groups acts on $\mod_\Lambda (\bd)$ by conjugations
\[
(g \cdot M)_\gamma = g_{t \gamma} M_\gamma g_{s \gamma}^{-1}, \,
\gamma \in \Delta_1,
\]
for $g \in \GL (\bd)$ and $M \in \mod_\Lambda (\bd)$. The orbits
with respect to this action correspond bijectively to the
isomorphism classes of $\Lambda$-modules of dimension vector
$\bd$.

\subsection{} \label{subsectZ1}
We present now a construction investigated in~\cite{Bon2}*{2.1} by
Bongartz. Fix $\bd', \bd'' \in \bbN^{\Delta_0}$, $M' \in
\mod_\Lambda (\bd')$ and $M'' \in \mod_\Lambda (\bd'')$. For $Z
\in \bbA (\bd', \bd'')$ and a path $\sigma = \gamma_1 \cdots
\gamma_m$ of positive length, let
\[
Z_\sigma = \sum_{i \in [1, m]} M_{\gamma_1}' \cdots M_{\gamma_{i -
1}}' Z_{\gamma_i} M_{\gamma_{i + 1}}'' \cdots M_{\gamma_m}''.
\]
If $\rho = \sum_{i \in [1, n]} \lambda_i \sigma_i$ is a relation
in $\Delta$, then $Z_\rho = \sum_{i \in [1, n]} \lambda_i
Z_{\sigma_i}$. We define $Z (M'', M')$ as the set of all $Z \in
\bbA (\bd', \bd'')$ such that $Z_\rho = 0$ for all $\rho \in R$.
For $Z \in Z (M'', M')$, let $M \in \bbA (\bd' + \bd'', \bd' +
\bd'')$ be given by
\[
M_\gamma =
\begin{bmatrix}
M_\gamma' & Z_\gamma \\ 0 & M_\gamma''
\end{bmatrix}, \, \gamma \in \Delta_1.
\]
Then $M \in \mod_\Lambda (\bd' + \bd'')$ and we have a short exact
sequence
\[
0 \to M' \xrightarrow{f} M \xrightarrow{g} M'' \to 0,
\]
with the maps $f$ and $g$ given by the canonical injections
$k^{d_x'} \to k^{d_x' + d_x''}$, $x \in \Delta_0$, and the
canonical surjections $k^{d_x' + d_x''} \to k^{d_x''}$, $x \in
\Delta_0$, respectively. On the other hand, for every short exact
sequence $\varepsilon$ of the form
\[
0 \to M' \to M \to M'' \to 0,
\]
there exists a (non-unique) element of $Z (M'', M')$ such that the
corresponding short exact sequence is isomorphic to $\varepsilon$.
More precisely, the map $Z (M'', M') \to \Ext_\Lambda^1 (M'', M')$
described above is a surjective linear map. The kernel of this map
consists of $Z \in Z (M'', M')$ such that the corresponding
sequence splits, i.e., there exists $h \in \bbV (\bd', \bd'') =
\prod_{x \in \Delta_0} \bbM (d_x', d_x'')$ such that $Z_\gamma =
M_\gamma' h_{s \gamma} - h_{t \gamma} M_\gamma''$ for all $\gamma
\in \Delta_1$. Consequently,
\[
\dim_k Z (M'', M') = [M'', M']^1 - [M'', M'] + \sum_{x \in
\Delta_0} d_x' d_x''.
\]

\subsection{} \label{subsectZ}
Let $\bd \in \bbN^{\Delta_0}$ and $M \in \mod_\Lambda (\bd)$.
There is a natural inclusion of the tangent space $T_M
\mod_\Lambda (\bd)$ to $\mod_\Lambda (\bd)$ at $M$ into $Z (M, M)$
(see~\cite{Kr1}*{(2.7)}). If $[M, M]^2 = 0$, then this map is an
isomorphism. Indeed, we have a sequence of inequalities, which
implies the claim:
\begin{align*}
a (\bd) & = \bdim \GL (\bd) - \langle \bd, \bd \rangle = \bdim \GL
(\bd) - [M, M] + [M, M]^1
\\ %
& = \dim_k Z (M, M) \geq \dim_k T_M \mod_\Lambda (\bd) \geq \dim_M
\mod_\Lambda (\bd) \geq a (\bd),
\end{align*}
where $\dim_M \mod_\Lambda (\bd)$ denotes the dimension of
$\mod_\Lambda (\bd)$ at $M$, i.e., the maximum of the dimensions
of the irreducible components of $\mod_\Lambda (\bd)$ passing
through $M$. It also follows from the above calculations that if
$[M, M]^2 = 0$, then $\dim_M \mod_\Lambda (\bd) = a (\bd)$ and $M$
is a nonsingular point of $\mod_\Lambda (\bd)$ (see also~\cite{Ge}
for a general proof of the last assertion).

Using similar inequalities and the fact, that $Z (M, M)$ is the
tangent space to the corresponding (not necessarily reduced)
scheme (see~\cite{Vo}), one proves the following fact.

\begin{prop} \label{propcompinter}
If $[M, M]^2$ vanishes generically on $\mod_\Lambda (\bd)$, then
the variety $\mod_\Lambda (\bd)$ is a complete intersection of
dimension $a (\bd)$. Moreover, in the above situation $M \in
\mod_\Lambda (\bd)$ is nonsingular if and only if $[M, M]^2 = 0$.
\end{prop}

\begin{proof}
This is just a more general formulation of the fact proved
in~\cite{BobSk2}*{Section~1}.
\end{proof}

\subsection{}
Let $\bd'$ and $\bd''$ be dimension vectors. Put $\bd = \bd' +
\bd''$. Let $C'$ be a constructible irreducible $\GL
(\bd')$-invariant subset of $\mod_\Lambda (\bd')$ and let $C''$ be
a constructible irreducible $\GL (\bd'')$-invariant subset of
$\mod_\Lambda (\bd'')$. Let
\begin{align*}
\hom (C', C'') & = \min \{ [M', M''] \mid M' \in C', \, M'' \in
C'' \},
\\ %
\ext^1 (C', C'') & = \min \{ [M', M'']^1 \mid M' \in C', \, M''
\in C'' \},
\\ %
\intertext{and} %
\ext^2 (C', C'') & = \min \{ [M', M'']^2 \mid M' \in C', \, M''
\in C'' \}.
\end{align*}
Recall from~\cite{CBSc}*{Lemma~4.3} that the functions
\begin{align*}
C' \times C'' \ni (M', M'') & \mapsto [M', M''] \in \bbZ,
\\ %
C' \times C'' \ni (M', M'') & \mapsto [M', M'']^1 \in \bbZ
\end{align*}
are upper semicontinuous. Moreover, in our case
\[
[M', M'']^2 = \dim_k Z (M', M'') - \sum_{\gamma \in \Delta_1} d_{s
\gamma}' d_{t \gamma}'' + \sum_{\rho \in R} d_{s \rho}' d_{t
\rho}'',
\]
hence the function
\[
C' \times C'' \ni (M', M'') \mapsto [M', M'']^2 \in \bbZ
\]
is also upper semicontinuous (using standard projective
resolutions one may prove this fact in a more general setting). In
particular, the sets
\begin{gather*}
\{ (M', M'') \in C' \times C'' \mid [M', M''] = \hom (C', C'') \},
\\ %
\{ (M', M'') \in C' \times C'' \mid [M', M'']^1 = \ext^1 (C', C'')
\},
\\ %
\{ (M', M'') \in C' \times C'' \mid [M', M'']^2 = \ext^2 (C', C'')
\}
\end{gather*}
are open subsets of $C' \times C''$.

We define $C' \oplus C''$ to be the set of all $M \in \mod_\Lambda
(\bd)$ which are isomorphic to a module of the form $M' \oplus
M''$ for $M' \in C'$ and $M'' \in C''$. We have the following
formula for the dimension of $C' \oplus C''$.

\begin{lemm}
If $C'$ and $C''$ are as above, then $C' \oplus C''$ is a
constructible $\GL (\bd)$-invariant irreducible subset of
$\mod_\Lambda (\bd)$ of dimension
\begin{multline*}
\dim C' + \dim C'' + \dim \GL (\bd)
\\ %
- \dim \GL (\bd') - \dim \GL (\bd'') - \hom  (C', C'') - \hom
(C'', C').
\end{multline*}
\end{lemm}

\begin{proof}
The claim follows by considering the map
\[
\GL (\bd) \times C' \times C'' \ni (g, M', M'') \mapsto g \cdot
(M' \oplus M'') \in \mod_\Lambda (\bd)
\]
(compare for example~\cite{CBSc}*{Section~1}).
\end{proof}

A special case of the above lemma, which is really of interest for
us, is the following.

\begin{coro} \label{corooplusdim}
Let $C'$ and $C''$ be as above. If $\dim C' = a (\bd') - c_1$,
$\dim C'' = a (\bd'') - c_2$, $\hom (C', C'') = \langle \bd',
\bd'' \rangle$ and $\hom (C'', C') = 0$, then $C' \oplus C''$ is a
constructible irreducible subset of $\mod_\Lambda (\bd)$ of
dimension
\[
a (\bd) + \langle \bd'', \bd' \rangle - (c_1 + c_2).
\]
\end{coro}

\begin{proof}
Direct calculations.
\end{proof}

\subsection{} \label{subsectext}
Let $C'$ and $C''$ be as above.  By $\calE (C', C'')$ we mean the
set of all $M \in \mod_\Lambda (\bd' + \bd'')$ for which there
exists an exact sequence
\[
0 \to M'' \to M \to M' \to 0
\]
with $M' \in C'$ and $M'' \in C''$. It follows
from~\cite{CBSc}*{Theorem~1.3(i)}, that if $C'$ and $C''$ are
closed subsets of $\mod_\Lambda (\bd')$ and $\mod_\Lambda (\bd'')$
respectively, then $\calE (C', C'')$ is a closed subset of
$\mod_\Lambda (\bd' + \bd'')$.

\subsection{}
Let $\bd$ be a dimension vector. Let $\mod_\Lambda^P$ be the full
subcategory of $\Lambda$-modules of projective dimension at most
$1$. Barot and Schr\"oer proved in~\cite{BaSc}*{Proposition~3.1}
that if $\mod_\Lambda^P (\bd)$ is nonempty, then it is an
irreducible open subset of $\mod_\Lambda (\bd)$ of dimension $a
(\bd)$. Here, for a subcategory $\calX$ of $\mod_\Lambda$ and $\bd
\in \bbN^{\Delta_0}$, we denote by $\calX (\bd)$ the set of all $M
\in \mod_\Lambda (\bd)$ such that $M \in \calX$. Dually,
$\mod_\Lambda^I (\bd)$ (if nonempty) is an irreducible open subset
of $\mod_\Lambda (\bd)$ of dimension $a (\bd)$, where
$\mod_\Lambda^I$ is the full subcategory of $\Lambda$-modules of
injective dimension at most $1$.

\subsection{}
From now till the end of the section we assume that $\Lambda$ is a
fixed canonical algebra. The first observation is the following.

\begin{lemm}
If $\bd \in \bP + \bR$, then $(\calP \vee \calR) (\bd)$ is an open
subset of $\mod_\Lambda^P (\bd)$. In particular, $\dim (\calP \vee
\calR) (\bd) = a (\bd)$.
\end{lemm}

By duality, if $\bd \in \bR + \bQ$, then $(\calR \vee \calQ)
(\bd)$ is an open subset of $\mod_\Lambda^I (\bd)$ of dimension $a
(\bd)$. As a consequence it also follows that, if $\bd \in \calR$,
then $\calR (\bd) = (\calP \vee \calR) (\bd) \cap (\calR \vee
\calQ) (\bd)$ is an irreducible open subset of $\mod_\Lambda
(\bd)$ of dimension $a (\bd)$ (see also~\cite{DoLe2}*{Section~4}
for another explanation of the last fact).

\begin{proof}
We already know that $(\calP \vee \calR) (\bd)$ is contained in
$\mod_\Lambda^P (\bd)$, thus it only remains to show that it is
open. But $M \in (\calP \vee \calR) (\bd)$ if and only if there
exists $X \in \calR$ of dimension vector $\bh$ such that
$\Hom_\Lambda (X, M) = 0$, hence the claim follows.
\end{proof}

\subsection{}
The proof of an analogous fact for $\bd \in \bP$ is more involved.

\begin{lemm}
If $\bd \in \bP$, then $\calP (\bd)$ is an open subset of
$\mod_\Lambda^P (\bd)$. In particular, $\dim \calP (\bd) = a
(\bd)$.
\end{lemm}

By duality, if $\bd \in \bQ$, then $\calQ (\bd)$ is an open subset
of $\mod_\Lambda^I (\bd)$ of dimension $a (\bd)$.

\begin{proof}
Again we only have to show that $\calP (\bd)$ is an open subset of
$\mod_\Lambda (\bd)$, hence also of $\mod_\Lambda^P (\bd)$. We
prove that $\calP (\bd)$ is the complement of the sum
\[
\bigcup_{\substack{\bd' \in \bP, \, \bd'' \in \bR + \bQ \\
\bd' + \bd'' = \bd, \, \bd'' \neq 0}} \calE (\mod_\Lambda (\bd'),
\mod_\Lambda (\bd'')).
\]
Since this is a finite sum of sets which are closed
by~\ref{subsectext}, it will imply the lemma. In order to show the
above claim, take $M \not \in \calP (\bd)$. Then $M = M' \oplus
M''$ for some $M' \in \calP$ and $M'' \in \calR \vee \calQ$, $M''
\neq 0$, and obviously
\[
M \in \calE (\mod_\Lambda (\bdim M'), \mod_\Lambda (\bdim M'')).
\]
Assume now that $M \in \calE (\mod_\Lambda (\bd'), \mod_\Lambda
(\bd''))$ for $\bd'$ and $\bd''$ as above. Let
\[
0 \to M'' \to M \to M' \to 0
\]
be a short exact sequence with $M' \in \mod_\Lambda (\bd')$ and
$M'' \in \mod_\Lambda (\bd'')$. Since $\langle \bdim M'', \bh
\rangle \leq 0$, it follows that $M''$ has a nonzero direct
summand $N'' \in \calR \vee \calQ$. Since $\Hom_\Lambda (N'', N) =
0$ for all $N \in \calP$, we get $M \not \in \calP$, and we are
done.
\end{proof}

\makeatletter
\def\@secnumfont{\mdseries} 
\makeatother

\section{Proofs of Theorems~\ref{theoirr} and~\ref{theocomp}} \label{sectproof}

\subsection{}
We prove Theorems~\ref{theoirr} and~\ref{theocomp} in a more
general setting. Let $\Lambda$ be the path algebra of a bound
quiver $(\Delta, R)$ of global dimension at most $2$. We also
assume that we are given two full subcategories $\calX$ and
$\calY$ of $\mod_\Lambda$ having the following properties:
\begin{enumerate}

\item
$\calX$ and $\calY$ are closed under forming direct sums and
taking direct summands,

\item
$\calX \vee \calY =  \mod_\Lambda$,

\item
$\pd_\Lambda M \leq 1$ for $M \in \calX$ and $\id_\Lambda N \leq
1$ for $N \in \calY$,

\item
$[N, M] = 0$ and $[M, N]^1 = 0$ for $N \in \calY$ and $M \in
\calX$,

\item
if $\bd \in \bbN^{\Delta_0}$, then $\calX (\bd)$ and $\calY (\bd)$
are open subsets of $\mod_\Lambda (\bd)$.

\end{enumerate}
Observe that canonical algebras fit into the above setting with
$\calX = \calP$ and $\calY = \calR \vee \calQ$ (or $\calX = \calP
\vee \calR$ and $\calY = \calQ$).

Let $\bX$ and $\bY$ denote the sets of the dimension vectors of
the modules belonging to $\calX$ and $\calY$, respectively. It
follows from the above conditions that if $\bd \in \bX$ ($\bd \in
\bY$, respectively), then $\calX (\bd)$ ($\calY (\bd)$) is an
irreducible open subset of $\mod_\Lambda (\bd)$ of dimension $a
(\bd)$. In particular, if $\bd' \in \bX$ and $\bd'' \in \bY$, then
$\calX (\bd') \oplus \calY (\bd'')$ is an irreducible
constructible subset (in fact,
using~\cite{CBSc}*{Theorem~1.3(iii)} one can even show that this
set is locally closed) of $\mod_\Lambda (\bd' + \bd'')$ of
dimension $a (\bd' + \bd'') + \langle \bd'', \bd' \rangle$ (apply
Corollary~\ref{corooplusdim}). Consequently, for $\bd \in
\bbN^{\Delta_0}$, $\mod_\Lambda (\bd)$ is a finite disjoint union
\[
\bigcup_{\substack{\bd' \in \bX, \, \bd'' \in \bY \\
\bd' + \bd'' = \bd}} \calX (\bd') \oplus \calY (\bd'')
\]
of irreducible constructible subsets of dimensions $a (\bd) +
\langle \bd'', \bd' \rangle$, respectively. In particular, this
implies that
\[
\dim \mod_\Lambda (\bd) = a (\bd) + \max \{ \langle \bd'', \bd'
\rangle \mid \bd' \in \bX, \, \bd'' \in \bY, \, \bd' + \bd'' = \bd
\}.
\]
Consequently, $\dim \mod_\Lambda (\bd) = a (\bd)$ if and only if
$\langle \bd'', \bd' \rangle \leq 0$ for all $\bd'$ and $\bd''$ as
above (recall that obviously $\dim \mod_\Lambda (\bd) \geq a
(\bd)$).

\subsection{}
As a first step in proving Theorem~\ref{theocomp} we prove the
following lemma.

\begin{lemm}
Let $\bd \in \bbN^{\Delta_0}$. If $\langle \bd'', \bd' \rangle
\leq 0$ for all $\bd' \in \bX$ and $\bd'' \in \bY$ such that $\bd
= \bd' + \bd''$, then $\mod_\Lambda (\bd)$ is a complete
intersection.
\end{lemm}

\begin{proof}
According to Proposition~\ref{propcompinter}, in order to prove
that $\mod_\Lambda (\bd)$ is a complete intersection, it is enough
to prove that $[M, M]^2$ vanishes generically on $\mod_\Lambda
(\bd)$. Note that every irreducible component of $\mod_\Lambda
(\bd)$ is the closure of the set of the form $\calX (\bd') \oplus
\calY (\bd'')$ for some $\bd' \in \bX$ and $\bd'' \in \bY$, such
that $\langle \bd'', \bd' \rangle = 0$. It is well-known that if
the closure of $\calX (\bd') \oplus \calY (\bd'')$ is an
irreducible component of $\mod_\Lambda (\bd)$, then $\ext^1 (\calY
(\bd''), \calX (\bd')) = 0$ (see for
example~\cite{CBSc}*{Theorem~1.2}). Since obviously $\hom (\calY
(\bd''), \calX (\bd')) = 0$, we get $\ext^2 (\calY (\bd''), \calX
(\bd')) = 0$ and the claim follows, because $\pd_\Lambda M' \leq
1$ for $M' \in \calX$ and $\id_\Lambda M'' \leq 1$ for $M'' \in
\calY $.
\end{proof}

\subsection{}
In order to reverse the above implication and finish the proof of
Theorem~\ref{theocomp} we need an additional assumption.

\begin{prop} \label{criterioncomplete}
Let $\bd$ be the dimension vector of a $\Lambda$-module of
projective or injective dimension at most $1$. Then $\mod_\Lambda
(\bd)$ is a complete intersection if and only if $\langle \bd'',
\bd' \rangle \leq 0$ for all $\bd' \in \bX$ and $\bd'' \in \bY$
such that $\bd = \bd' + \bd''$.
\end{prop}

\begin{proof}
We only have to prove that if $\mod_\Lambda (\bd)$ is a complete
intersection, then $\dim \mod_\Lambda (\bd) = a (\bd)$, but this
follows since $\mod_\Lambda (\bd)$ has an irreducible component of
dimension $a (\bd)$ (the closure of $\mod_\Lambda^P (\bd)$ or
$\mod_\Lambda^I (\bd)$) and complete intersections are
equidimensional.
\end{proof}

\subsection{}
We divide the proof of an analogous criterion for irreducibility
also in two steps.

\begin{lemm}
Let $\bd$ be the dimension vector. If $\langle \bd'', \bd' \rangle
\leq 0$ for all $\bd' \in \bX$ and $\bd'' \in \bY$ such that $\bd
= \bd' + \bd''$, and equality holds for exactly one pair $(\bd',
\bd'')$, then $\mod_\Lambda (\bd)$ is irreducible.
\end{lemm}

\begin{proof}
Let $\hat{\bd}' \in \bX$ and $\hat{\bd}'' \in \bY$ be such that
$\hat{\bd}' + \hat{\bd}'' = \bd$ and $\langle \hat{\bd}'',
\hat{\bd}' \rangle = 0$. Then $\calX (\hat{\bd}') \oplus \calY
(\hat{\bd}'')$ is an irreducible constructible subset of
$\mod_\Lambda (\bd)$ of dimension $a (\bd)$, and the remaining
sets $\calX (\bd') \oplus \calY (\bd'')$ have dimensions smaller
that $a (\bd)$. Since every irreducible component of $\mod_\Lambda
(\bd)$ has dimension at least $a (\bd)$, $\mod_\Lambda (\bd)$ is
the closure of $\calX (\hat{\bd}') \oplus \calY (\hat{\bd}'')$,
hence irreducible.
\end{proof}

\subsection{}
We may again reverse the above implication if we assume the
existence of a $\Lambda$-module of dimension vector $\bd$ and
projective or injective dimension at most $1$.

\begin{prop} \label{criterion}
Let $\bd$ be the dimension vector of a $\Lambda$-module of
projective or injective dimension at most $1$. Then $\mod_\Lambda
(\bd)$ is irreducible if and only if $\langle \bd'', \bd' \rangle
\leq 0$ for all $\bd' \in \bX$ and $\bd'' \in \bY$ such that $\bd
= \bd' + \bd''$, and equality holds for exactly one pair $(\bd',
\bd'')$.
\end{prop}

\begin{proof}
We only have to prove that if $\mod_\Lambda (\bd)$ is irreducible,
then the above condition is satisfied. Without loss of generality
we may assume that $\bd$ is the dimension vector of a
$\Lambda$-module of projective dimension at most $1$. Then we know
that the closure of $\mod_\Lambda^P (\bd)$ is an irreducible
component of $\mod_\Lambda (\bd)$ of dimension $a (\bd)$, thus
$\dim \mod_\Lambda (\bd) = a (\bd)$. In particular, $\langle
\bd'', \bd' \rangle \leq 0$ for all $\bd'$ and $\bd''$. Moreover,
the irreducible components of $\mod_\Lambda (\bd)$ are precisely
the closures of the sets $\calX (\bd') \oplus \calY (\bd'')$ with
$\dim (\calX (\bd') \oplus \calY (\bd'')) = a (\bd)$, i.e.,
$\langle \bd'', \bd' \rangle  = 0$. Thus irreducibility of
$\mod_\Lambda (\bd)$ implies that the equality holds for exactly
one pair $(\bd', \bd'')$.
\end{proof}

\subsection{}
As a consequence of the above propositions we also obtain some
information about the maximal $\GL (\bd)$-orbits in the above
situations. Namely, if $\bd$ is a dimension vector such that $\dim
\mod_\Lambda (\bd) = a (\bd)$, then every maximal $\GL
(\bd)$-orbit consists of points which are nonsingular in
$\mod_\Lambda (\bd)$. Indeed, let the $\GL (\bd)$-orbit of a
$\Lambda$-module $M$ be maximal and write $M = M' \oplus M''$ for
$M' \in \calX$ and $M'' \in \calY$. We know that $\langle \bdim
M'', \bdim M' \rangle \leq 0$. Obviously, $[M'', M'] = 0$.
Moreover, the maximality of the $\GL (\bd)$-orbit of $M$ implies
that $[M'', M']^1 = 0$ (see for example~\cite{Bon4}*{Lemma~1.1}).
Consequently, $[M'', M']^2 = 0$. Since $\pd_\Lambda M' \leq 1$ and
$\id_\Lambda M'' \leq 1$, this implies that $[M, M]^2 = 0$, which
finishes the proof according to~\ref{subsectZ}.

Assume now in addition that $\Lambda$ is a canonical algebra, $\bd
\in \bR$ and $\mod_\Lambda (\bd)$ is irreducible. Then it follows
from Proposition~\ref{criterion}, that $\langle \bd'', \bd'
\rangle < 0$ for all $\bd' \in \bP$ and $\bd'' \in \bR + \bQ$ with
$\bd' \neq 0$. Consequently, using again~\cite{Bon4}*{Lemma~1.1}
we obtain that if the $\GL (\bd)$-orbit of $M$ is maximal, then $M
\in \calR \vee \calQ$. Since $\langle \bh, \bd \rangle = 0$,
$(\calR \vee \calQ) (\bd) = \calR (\bd)$ and $M \in \calR$. With
methods analogous to those used in the proofs
of~\cite{Ri2}*{Theorem~3.5} and~\cite{BobSk2}*{Proposition~5}, one
can give a precise description of the maximal $\GL (\bd)$-orbits.
It is essentially identical to that given
in~\cite{BobSk2}*{Proposition~5}, but since it is lengthy and
requires introducing an appropriate language, we will not present
it here.

\subsection{} \label{lemmext2}
We give now the proof Theorem~\ref{theoirr}. The crucial
observation, whose proof is based on ideas of the proof
of~\cite{Re}*{Proposition~2.5}, is the following.

\begin{lemm}
Let $\bd' \in \bX$ and $\bd'' \in \bY$. If $\dim \mod_\Lambda
(\bd' + \bd'') = a (\bd' + \bd'')$, then $\ext^1 (\calY (\bd''),
\calX (\bd')) = - \langle \bd'', \bd' \rangle$.
\end{lemm}

\begin{proof}
Let $\bd = \bd' + \bd''$, $C' = \calX (\bd')$ and $C'' = \calY
(\bd'')$. Obviously, $\ext^1 (C'', C') \geq - \langle \bd'', \bd'
\rangle$, thus we only have to show that $\ext^1 (C'', C') \leq -
\langle \bd'', \bd' \rangle$.

Recall that $\calU = \{ (M', M'') \in C' \times C'' \mid [M'',
M']^1 = \ext^1 (C'', C') \}$ is an open subset of $C' \times C''$.
Consequently, according to~\cite{CBSc}*{Theorem~1.3(iii)}, the
subset $\calV$ of all $M \in \mod_\Lambda (\bd)$ such that there
exists a short exact sequence
\[
0 \to M' \to M \to M'' \to 0
\]
with $(M', M'') \in \calU$ is an open subset of $\mod_\Lambda
(\bd)$. In particular, $\calV$ contains nonsingular points of
$\mod_\Lambda (\bd)$.

Let $\calZ = \{ (M', M'', Z) \mid (M', M'') \in \calU, \, Z \in Z
(M'', M') \}$. It follows from~\cite{Bon2}*{Lemma~1} that the
canonical projection $\calZ \to \calU$ is a subbundle of the
trivial vector bundle $\calU \times \bbA (\bd', \bd'') \to \calU$.
In particular, $\calZ$ is smooth, since $\calU$ is an open subset
of $\mod_\Lambda (\bd') \times \mod_\Lambda (\bd'')$ consisting of
nonsingular points. We describe now the tangent space $T_{(M',
M'', Z)} \calZ$ for $(M', M'', Z) \in \calZ$ more precisely.
Obviously, it is a subspace of
\[
\bbA (\bd', \bd') \times \bbA (\bd'', \bd'') \times \bbA (\bd',
\bd'')
\]
of dimension
\[
a (\bd') + a (\bd'') + \ext^1 (C'', C') + \sum_{x \in \Delta_0}
d_x' d_x''.
\]
Moreover, if $(Z', Z'', Y) \in T_{(M', M'', Z)} \calZ$, then $Z'
\in Z (M', M')$, $Z'' \in Z (M'', M'')$ and $Y_\rho = r_\rho$ for
all $\rho \in R$, where $Y_\rho$ is defined as in
Section~\ref{subsectZ} and we define $r_\rho$ by the standard
extension to relations of the following definition for paths: if
$\sigma = \gamma_1 \cdots \gamma_m$ is a path of length at least
$2$, then
\begin{multline*}
r_\sigma = \sum_{i < j \in [1, m]} M_{\gamma_1}' \cdots
M_{\gamma_{i - 1}}' Z_{\gamma_i}' M_{\gamma_{i + 1}}' \cdots
M_{\gamma_{j - 1}}' Z_{\gamma_j} M_{\gamma_{j + 1}}'' \cdots
M_{\gamma_m}'' +
\\ %
\sum_{i < j \in [1, m]} M_{\gamma_1}' \cdots M_{\gamma_{i - 1}}'
Z_{\gamma_i} M_{\gamma_{i + 1}}'' \cdots M_{\gamma_{j - 1}}''
Z_{\gamma_j}'' M_{\gamma_{j + 1}}'' \cdots M_{\gamma_m}''.
\end{multline*}
Since $\dim_k Z (M', M') = a (\bd')$, $\dim_k Z (M'', M'') = a
(\bd'')$ and the solution set of the homogeneous system $Y_\rho =
0$, $\rho \in R$, is $Z (M'', M')$, it follows by comparing
dimensions that $T_{(M', M'', Z)} \calZ$ is the set of all $(Z',
Z'', Y)$ satisfying the above conditions.

Note that $T_{\Id} \GL (\bd) = \bbV (\bd, \bd)$
(see~\ref{subsectZ1} for a definition). With respect to the
canonical decomposition
\[
\bbV (\bd, \bd) =
\begin{bmatrix}
\bbV (\bd', \bd') & \bbV (\bd', \bd'') \\ \bbV (\bd'', \bd') &
\bbV (\bd'', \bd'')
\end{bmatrix},
\]
every element $X \in \bbV (\bd, \bd)$ can be written as
\[
X =
\begin{bmatrix}
X^{(1, 1)} & X^{(1, 2)} \\ X^{(2, 1)} & X^{(2, 2)}
\end{bmatrix}.
\]
On the other hand, if $M \in \calV$, then $T_M \calV = T_M
\mod_\Lambda (\bd) \subset \bbA (\bd, \bd)$, and the canonical
decomposition
\[
\bbA (\bd, \bd) =
\begin{bmatrix}
\bbA (\bd', \bd') & \bbA (\bd', \bd'') \\ \bbA (\bd'', \bd') &
\bbA (\bd'', \bd'')
\end{bmatrix}
\]
induces the analogous matrix presentation of the elements of $T_M
\calV$.

We have a surjective map $\Phi : \GL (\bd) \times \calZ \to \calV$
given by
\[
(g, M', M'', Z) \mapsto g \cdot M, \text{ where } M_\gamma =
\begin{bmatrix}
M_\gamma' & Z_\gamma \\ 0 & M_\gamma''
\end{bmatrix}, \, \gamma \in \Delta_1.
\]
For fixed $(M', M'', Z) \in \calZ$, we have the tangent map $F :
T_{\Id} \GL (\bd) \times T_{(M', M'', Z)} \calZ \to T_M \calV$
given by (here we apply the conventions about presenting the
elements of $\bbV (\bd, \bd)$ and $T_M \calV$ introduced above)
\begin{align*}
F (X, Z', Z'', Y)_\gamma^{(1, 1)} & = Z_\gamma' + X_{t
\gamma}^{(1, 1)} M_\gamma' - M_\gamma' X_{s \gamma}^{(1, 1)} -
Z_\gamma X_{s \gamma}^{(2, 1)},
\\ %
F (X, Z', Z'', Y)_\gamma^{(1, 2)} & = Y_\gamma + X_{t \gamma}^{(1,
1)} Z_\gamma + X_{t \gamma}^{(1, 2)} M_\gamma''
\\ %
& \qquad - M_\gamma' X_{s \gamma}^{(1, 2)} - Z_\gamma X_{s
\gamma}^{(2, 2)},
\\ %
F (X, Z', Z'', Y)_\gamma^{(2, 1)} & = X_{t \gamma}^{(2, 1)}
M_\gamma' - M_\gamma'' X_{s \gamma}^{(2, 1)},
\\ %
\intertext{and} %
F (X, Z', Z'', Y)_\gamma^{(2, 2)} & =  Z_\gamma'' + X_{t
\gamma}^{(2, 1)} Z_\gamma + X_{t \gamma}^{(2, 2)} M_\gamma'' -
M_\gamma'' X_{s \gamma}^{(2, 2)},
\end{align*}
for $\gamma \in \Delta_1$ (it follows by computations using block
matrices --- compare the corresponding calculations in the proof
of~\cite{Re}*{Proposition~2.5}).

Note that $X^{(2, 1)} \in \Hom_\Lambda (M', M'')$ for $(X, Z',
Z'', Y) \in \Ker F$. Thus the linear map $G : \Ker F \to
\Hom_\Lambda (M', M'')$ given by $G (X, Z', Z'', Y) = X^{(2, 1)}$
is well-defined. Moreover, $G$ is surjective. Indeed, for $f : M'
\to M''$ we define
\begin{align*}
X & =
\begin{bmatrix}
0 & 0 \\ f & 0
\end{bmatrix},
& Y_\gamma & = 0,
\\ %
Z_\gamma' & = Z_\gamma f_{s \gamma}, & Z_\gamma'' & = - f_{t
\gamma} Z_\gamma,
\end{align*}
where $\gamma \in \Delta_1$. One checks that $(Z', Z'', Y) \in
T_{(M', M'', Z)} \calZ$. Moreover, $F (X, Z', Z'', Y) = 0$ and $G
(X, Z', Z'', Y) = f$.

Let $\frakp = \{ X \in \bbV (\bd) \mid X^{(2, 1)} = 0 \}$. Define
$H : \frakp \to T_{\Id} \GL (\bd) \times T_{(M', M'', Z)} \calZ$
by $H (X) = (X, Z', Z'', Y)$, where
\begin{align*}
Z_\gamma' & =  M_\gamma' X_{s \gamma}^{(1, 1)} - X_{t \gamma}^{(1,
1)} M_\gamma',
\\ %
Z_\gamma'' & = M_\gamma'' X_{s \gamma}^{(2, 2)} - X_{t
\gamma}^{(2, 2)} M_\gamma'',
\\ %
\intertext{and} %
Y_\gamma & = M_\gamma' X_{s \gamma}^{(1, 2)} + Z_\gamma X_{s
\gamma}^{(2, 2)} - X_{t \gamma}^{(1, 1)} Z_\gamma - X_{t
\gamma}^{(1, 2)} M_\gamma'',
\end{align*}
for $\gamma \in \Delta_1$. Using the description of $T_{(M', M'',
Z)} \calZ$ one checks that $H$ is well-defined. Obviously, $H$ is
injective. Moreover, by direct calculations one checks that $\Im H
= \Ker G$. Consequently, we get
\begin{multline*}
\dim_k \Ker F = \dim_k \frakp + \dim_k \Hom_\Lambda (M', M'')
\\ %
= \dim \GL (\bd') + \dim \GL (\bd'') + \sum_{x \in \Delta_0} d_x'
d_x'' + \langle \bd', \bd'' \rangle.
\end{multline*}

We may assume that $M$ is a nonsingular point of $\mod_\Lambda
(\bd)$. Then we have the following sequence of inequalities
\begin{align*}
a (\bd) & = \dim_M \mod_\Lambda (\bd) = \dim_k T_M \mod_\Lambda
(\bd) = \dim_k T_M \calV
\\ %
& \geq \dim_k \Im F = \dim_k \bbV (\bd) + \dim_k T_{(M', M'', Z)}
\calZ - \dim_k \Ker F
\\ %
& = \dim \GL (\bd) + a (\bd') + a (\bd'') + \ext^1 (C'', C') +
\sum_{x \in \Delta_0} d_x' d_x'' -
\\ %
& \qquad - \dim \GL (\bd') - \dim \GL (\bd'') - \sum_{x \in
\Delta_0} d_x' d_x'' - \langle \bd', \bd'' \rangle
\\ %
& = a (\bd) + \langle \bd'', \bd' \rangle + \ext^1 (C'', C'),
\end{align*}
which implies $\ext^1 (C'', C') \leq - \langle \bd'', \bd'
\rangle$, hence finishes the proof.
\end{proof}

\subsection{}
Another useful observation is the following.

\begin{lemm}
Let $\bd$ be a dimension vector such that $\mod_\Lambda (\bd)$ is
irreducible, $\mod_\Lambda^P (\bd) \neq \varnothing$ and
$\mod_\Lambda^I (\bd) \neq \varnothing$. If $\bd' \in \bX$ and
$\bd'' \in \bY$ are such that $\bd' + \bd'' = \bd$ and $\langle
\bd'', \bd' \rangle = 0$, then the set of $M \in \calX (\bd')
\oplus \calY (\bd'')$ such that $[M, M]^2 \neq 0$ has dimension at
most $a (\bd) - 2$.
\end{lemm}

\begin{proof}
Our assumptions imply that $\mod_\Lambda (\bd)$ is the closure of
$\calX (\bd') \oplus \calY (\bd'')$, hence $(\calX (\bd') \oplus
\calY (\bd'')) \cap \mod_\Lambda^P (\bd) \neq \varnothing$ and
$(\calX (\bd') \oplus \calY (\bd'')) \cap \mod_\Lambda^I (\bd)
\neq \varnothing$. In particular, $\calX (\bd') \cap
\mod_\Lambda^I (\bd') \neq \varnothing$ and $\calY (\bd'') \cap
\mod_\Lambda^P (\bd'') \neq \varnothing$. Consequently, $\dim C'
\leq a (\bd') - 1$ and $\dim C'' \leq a (\bd'') - 1$, where $C' =
\{ M' \in \calX (\bd') \mid \id_\Lambda M' = 2 \}$ and $C'' = \{
M'' \in \calY (\bd'') \mid \pd_\Lambda M'' = 2 \}$. Since $\{ M
\in \calX (\bd') \oplus \calY (\bd'') \mid [M, M]^2 \neq 0 \}
\subset C' \oplus C''$, the claim follows from
Corollary~\ref{corooplusdim}.
\end{proof}

\subsection{} \label{propnorm}
The following fact implies Theorem~\ref{theoirr}.

\begin{prop}
Let $\bd$ be a dimension vector such that at least one of the
following conditions is satisfied:
\begin{enumerate}

\item
$\bd \in \bX$ or $\bd \in \bY$,

\item
$\mod_\Lambda^P (\bd) \neq \varnothing$ and $\mod_\Lambda^I (\bd)
\neq \varnothing$.

\end{enumerate}
Then $\mod_\Lambda (\bd)$ is normal if and only if it is
irreducible.
\end{prop}

\begin{proof}
We only have to prove that if $\mod_\Lambda (\bd)$ is irreducible,
then it is normal. First observe, that irreducibility of
$\mod_\Lambda (\bd)$ implies that either $\mod_\Lambda^P (\bd)$ or
$\mod_\Lambda^I (\bd)$ is a dense open subset of $\mod_\Lambda
(\bd)$. In particular, $[M, M]^2$ vanishes generically on
$\mod_\Lambda (\bd)$, hence $\mod_\Lambda (\bd)$ is a complete
intersection by Proposition~\ref{propcompinter}. Consequently,
according to the Serre's criterion (see for
example~\cite{Ha}*{Proposition~8.23}) in order to prove normality
we have to show that $\mod_\Lambda (\bd)$ is nonsingular in
codimension $1$. According to Proposition~\ref{propcompinter},
this will follow if we show that the set of $M \in \mod_\Lambda
(\bd)$ such that $[M, M]^2 \neq 0$ is of codimension at least $2$.

Recall that
\[
\mod_\Lambda (\bd) = \bigcup_{\substack{\bd' \in \bX, \, \bd'' \in
\bY \\ \bd' + \bd'' = \bd}} \calX (\bd') \oplus \calY (\bd'')
\]
is a presentation of $\mod_\Lambda (\bd)$ as a finite disjoint sum
of constructible sets of dimensions $a (\bd) + \langle \bd'', \bd'
\rangle$, respectively. Thus we have to show that the set of $M
\in \calX (\bd') \oplus \calY (\bd'')$ such that $[M, M]^2 \neq 0$
has dimension at most $a (\bd) - 2$ for all $\bd'$ and $\bd''$ as
above. Note that by Lemma~\ref{lemmext2} $\ext^1 (\calY (\bd''),
\calX (\bd')) = -\langle \bd'', \bd' \rangle$. Since obviously
$\hom (\calY (\bd''), \calX (\bd')) = 0$, we get $\ext^2 (\calY
(\bd''), \calX (\bd')) = 0$. This implies our claim if $\langle
\bd'', \bd' \rangle < 0$.

Assume now that $\bd' \in \bX$, $\bd'' \in \bY$, $\bd' + \bd'' =
\bd$ and $\langle \bd'', \bd' \rangle = 0$. If $\bd \in \bX$, then
Proposition~\ref{criterion} implies, that $\bd' = \bd$ and $\bd''
= 0$. Consequently, $[M, M]^2 = 0$ for all $M \in \calX (\bd')
\oplus \calY (\bd'') = \calX (\bd)$. A similar argument applies if
$\bd \in \bY$. If $\mod_\Lambda^P (\bd) \neq \varnothing$ and
$\mod_\Lambda^I (\bd) \neq \varnothing$, then we can use the
previous lemma.
\end{proof}

\makeatletter
\def\@secnumfont{\mdseries} 
\makeatother

\section{Inequalities} \label{sectineq}

Throughout this section, $\Lambda$ is a fixed canonical algebra of
type $\bm = (m_1, \ldots, m_n)$.

\subsection{}
Our aim in this section is to prove the following inequalities.

\begin{prop} \label{propineqnew}
Let $\bd \in \bR$ and $\bd' \in \bP$ be such that $\bd' \neq 0$
and $\bd - \bd' \in \bR + \bQ$.
\begin{enumerate}

\item
If $\sum_{i \in [1, n]} \frac{1}{m_i - 1} > 2 n - 5$, then
$\langle \bd - \bd', \bd' \rangle < 0$.

\item
If $\sum_{i \in [1, n]} \frac{1}{m_i - 1} = 2 n - 5$, then
$\langle \bd - \bd', \bd' \rangle \leq 0$. Moreover, if $p^{\bd} >
0$, then the above inequality is strict.

\item
If $\bm = (2, 2, 2, 2, 2)$ and $\bd$ is sincere, then $\langle \bd
- \bd', \bd' \rangle < 0$.

\end{enumerate}
\end{prop}

According to the results of the previous paragraph, the above
proposition implies the ``positive'' parts of
Theorems~\ref{theo1}, \ref{theo2} and \ref{theo3}.

\subsection{}
We start with a series of simple inequalities leading to our main
result. The elementary proof of the first inequality is left to
the reader.

\begin{lemm}
Let $d$ and $\delta_1$, \ldots, $\delta_m$, $m > 0$, be
nonnegative and such that $d = \sum_{i \in [1, m]} \delta_i$. Then
\[
\sum_{i < j \in [1, m]} \delta_i \delta_j \leq \ffrac{m - 1}{2 m}
d^2.
\]
Moreover, equality holds if and only if $\delta_i = \frac{d}{m}$
for all $i \in [1, m]$. \qed
\end{lemm}

\subsection{}
We will need the following variant of the above inequality.

\begin{lemm}
Let $d$ and $\delta_1$, \ldots, $\delta_m$, $m > 2$, be
nonnegative and such that $d = \sum_{i \in [1, m]} \delta_i$. Then
\[
\sum_{i < j \in [1, m]} \delta_i \delta_j \leq \ffrac{1}{4} (d +
d')^2 - \ffrac{m - 1}{2 (m - 2)} d'^2,
\]
where $d' = \sum_{i \in [2, m - 1]} \delta_i$. Moreover, equality
holds if and only if $\delta_1 = \delta_m = \frac{d - d'}{2}$ and
$\delta_i = \frac{d'}{m - 2}$ for all $i \in [2, m - 1]$.
\end{lemm}

\begin{proof}
It follows by direct calculations that
\[
\sum_{i < j \in [1, m]} \delta_i \delta_j = \delta_1 \delta_m + (d
- d') d' + \sum_{i < j \in [2, m - 1]} \delta_i \delta_j.
\]
Now the claim follows by applying the previous lemma to $\delta_1
\delta_m$ and $\sum_{i < j \in [2, m - 1]} \delta_i \delta_j$.
\end{proof}

\subsection{}
The next step is the following.

\begin{lemm}
Let $d$, $q$ and $\delta_1$, \ldots, $\delta_m$, $m \geq 2$, be
nonnegative and such that $d = \sum_{i \in [1, m]} \delta_i$ and
$\delta_1, \delta_m \geq q$. Then
\[
\sum_{i < j \in [1, m]} \delta_i \delta_j \leq
\begin{cases}
\frac{m - 1}{2 m} d^2 & m q \leq d,
\\ %
(d - q)^2 - \frac{m - 1}{2 (m - 2)} (d - 2 q)^2 & m q > d.
\end{cases}
\]
Moreover, in the first case equality holds if and only if
$\delta_i = \frac{d}{m}$ for all $i \in [1, m]$, and in the second
case equality holds if and only if $\delta_1 = \delta_m = q$ and
$\delta_i = \frac{d - 2 q}{m - 2}$ for all $i \in [2, m - 1]$.
\end{lemm}

Note that $m q > d$ may hold only for $m > 2$.

\begin{proof}
The claim for $m = 2$ is an easy exercise, hence we may assume
that $m > 2$ and apply the previous lemma. Since $d' = \sum_{i \in
[2, m - 1]} \delta_i$ varies from $0$ to $d - 2 q$, the maximal
value of $\frac{1}{4} d^2 + \frac{1}{2} d d' - \frac{m}{4 (m - 2)}
d'^2$ is obtained for $d' = \min (\frac{m - 2}{m} d, d - 2 q)$.
This immediately implies our claim.
\end{proof}

\subsection{}
The following inequality is what we really need.

\begin{lemm} \label{lemmineq}
Let $d$, $q$ and $\delta_1$, \ldots, $\delta_m$, $m \geq 2$, be
nonnegative and such that $d = \sum_{i \in [1, m]} \delta_i$ and
$\delta_1 \geq q$. Then
\[
- \delta_m q + \sum_{i < j \in [1, m]} \delta_i \delta_j \leq
\begin{cases}
- d q + \frac{m - 1}{2 m} (d + q)^2 & (m - 1) q \leq d,
\\ %
- d q + d^2 - \frac{m - 1}{2 (m - 2)} (d - q)^2 & (m - 1) q > d.
\end{cases}
\]
Moreover, in the first case equality holds if and only if
$\delta_i = \frac{d + q}{m}$ for all $i \in [1, m - 1]$ and
$\delta_m = \frac{d - (m - 1) q}{m}$, and in the second case
equality holds if and only if $\delta_1 = q$, $\delta_m = 0$ and
$\delta_i = \frac{d - q}{m - 2}$ for all $i \in [2, m - 1]$.
\end{lemm}

\begin{proof}
Let $\delta_i' = \delta_i$ for $i \in [1, m - 1]$ and $\delta_m' =
\delta_m + q$. Then
\[
- \delta_m q + \sum_{i < j \in [1, m]} \delta_i \delta_j = - d q +
\sum_{i < j \in [1, m]} \delta_i' \delta_j'
\]
and we may apply the previous lemma.
\end{proof}

\subsection{}
We will also need the following consequence of the previous
inequality.

\begin{coro} \label{coroineq}
Let $d$, $q$ and $\delta_1$, \ldots, $\delta_m$, $m \geq 2$, be
nonnegative and such that $d = \sum_{i \in [1, m]} \delta_i$ and
$\delta_1 \geq q$. Then
\[
- \delta_m q + \sum_{i < j \in [1, m]} \delta_i \delta_j \leq
\ffrac{1}{2} d^2 - \ffrac{1}{2} q^2.
\]
Moreover, if equality holds then $q = d$.
\end{coro}

\begin{proof}
If $q = 0$ then the claim is obvious from the previous lemma, thus
we assume that $q > 0$. We first consider the case $(m - 1) q \leq
d$. In this case $\frac{1}{m} \geq \frac{q}{d + q}$. Using once
more the previous lemma one easily gets that
\[
(\ffrac{1}{2} d^2 - \ffrac{1}{2} q^2) - \bigl( -\delta_m q +
\sum_{i < j \in [1, m]} \delta_i \delta_j \bigr) \geq \ffrac{1}{2
m} (d + q)^2 - q^2 \geq \ffrac{1}{2} q (d - q),
\]
hence the claim follows. On the other hand, if $(m - 1)q > d$,
then
\[
(\ffrac{1}{2} d^2 - \ffrac{1}{2} q^2) - \bigl( -\delta_m q +
\sum_{i < j \in [1, m]} \delta_i \delta_j \bigr) \geq \ffrac{1}{2
(m - 2)} (d - q)^2,
\]
which finishes the proof.
\end{proof}

\subsection{} \label{lemmfunf}
For a fixed positive $d$ and integers $m_1, m_2, m_3 \geq 2$, let
$f$ be the function defined on the set of all $4$-tuples $(p, p_1,
p_2, p_3)$ of nonnegative real numbers such that $p + p_1 + p_2 +
p_3 = d$ by $f (p, p_1, p_2, p_3) = g_{m_1} (p_1) + g_{m_2} (p_2)
+ g_{m_3} (p_3)$, where
\[
g_m (q)  =
\begin{cases}
\frac{m - 1}{2 m} (d + q)^2 & (m - 1) q \leq d,
\\ %
d^2 - \frac{m - 1}{2 (m - 2)} (d - q)^2 & (m - 1) q > d.
\end{cases}
\]
Our next aim is to prove the following.

\begin{lemm}
If $\frac{1}{m_1 - 1} + \frac{1}{m_2 - 1} + \frac{1}{m_3 - 1} \geq
1$, then $f (p, p_1, p_2, p_3) \leq  2 d^2$. Moreover, equality
holds if and only if $\frac{1}{m_1 - 1} + \frac{1}{m_2 - 1} +
\frac{1}{m_3 - 1} = 1$, $p = 0$ and there exists $i \in [1, 3]$
and $\lambda$ such that $\frac{m_i - 1}{m_i} d \leq \lambda \leq
d$, $p_i = \frac{m_i}{m_i - 1} \lambda - d$ and $p_j = d -
\frac{m_j - 2}{m_j - 1} \lambda$ for $j \neq i$.
\end{lemm}

\begin{proof}
If we substitute $p = \xi^2$ and $p_i = \xi_i^2$ for $i \in [1,
3]$, then we may replace $f$ by a function $F$ defined on the set
of all $4$-tuples $(\xi, \xi_1, \xi_2, \xi_3)$ of real numbers
such that $\xi^2 + \xi_1^2 + \xi_2^2 + \xi_3^2 = d$ by $F (\xi,
\xi_1, \xi_2, \xi_3) = G_{m_1} (\xi_1) + G_{m_2} (\xi_2) + G_{m_3}
(\xi_3)$, where
\[
G_m (\mu)  =
\begin{cases}
\frac{m - 1}{2 m} (d + \mu^2)^2 & (m - 1) \mu^2 \leq d,
\\ %
d^2 - \frac{m - 1}{2 (m - 2)} (d - \mu^2)^2 & (m - 1) \mu^2 > d.
\end{cases}
\]
By direct calculations one checks that $F$ is differentiable.
Since the set considered is compact, $F$ posses a maximum. Using
Lagrange's multipliers method we know that, if $F$ has a maximum
at $(\xi, \xi_1, \xi_2, \xi_3)$, then there exists $\lambda$ such
that $\lambda \xi = 0$ and $\xi_i H_{m_i} (\xi_i) = 0$ for all $i
\in [1, 3]$, where
\[
H_m (\mu) =
\begin{cases}
\frac{m - 1}{m} (d + \mu^2) - \lambda & (m - 1) \mu^2 \leq d,
\\ %
\frac{m - 1}{m - 2} (d - \mu^2) - \lambda & (m - 1) \mu^2 > d.
\end{cases}
\]
If $\xi \neq 0$, then $\lambda = 0$  and it follows that either
$\xi_i = 0$ or $\xi_i^2 = d$ for each $i$. Let $I$ be the set of
all $i$ such that $\xi_i^2 = d$. Then $d |I| + \xi = d$, hence $I
= \varnothing$ and $\xi^2 = d$. We have
\[
f (d, 0, 0, 0) = (\ffrac{3}{2} - \ffrac{1}{2 m_1} - \ffrac{1}{2
m_2} - \ffrac{1}{2 m_3}) d^2 < 2 d^2,
\]
thus we may assume that $\xi = 0$.

Let $I_0 = \{ i \mid \xi_i = 0 \}$, $I_1 = \{ i \mid 0 < \xi^2
\leq \frac{1}{m_i - 1} d \}$ and $I_2 = \{ i \mid \frac{1}{m_i -
1} d < \xi^2 \leq d \}$. Up to symmetry we have to consider the
following cases:
\begin{itemize}

\item
$I_0 = \{ 1, 2 \}$;

\item
$I_0 = \{ 1 \}$, $I_1 = \{ 2, 3 \}$;

\item
$I_0 = \{ 1 \}$, $I_1 = \{ 2 \}$, $I_2 = \{ 3 \}$;

\item
$I_0 = \{ 1 \}$, $I_2 = \{ 2, 3 \}$;

\item
$I_1 = \{ 1, 2, 3 \}$;

\item
$I_1 = \{ 1, 2 \}$, $I_2 = \{ 3 \}$;

\item
$I_1 = \{ 1 \}$, $I_2 = \{ 2, 3 \}$.

\end{itemize}
Note that our assumption $\xi = 0$ implies that $I_0 \neq \{ 1, 2,
3 \}$. On the other hand $\frac{1}{m_1 - 1} + \frac{1}{m_2 - 1} +
\frac{1}{m_3 - 1} \geq 1$ implies that $I_2 \neq \{ 1, 2, 3 \}$.
For future use, let $\delta = \frac{1}{m_1 - 1} + \frac{1}{m_2 -
1} + \frac{1}{m_3 - 1}$ and $\gamma = \frac{1}{m_2 - 1} +
\frac{1}{m_3 - 1}$.

We start with the case $I_0 = \{ 1, 2 \}$, thus $p_1 = p_2 = 0$.
Then obviously $p_3 = d$ and
\[
f (p, p_1, p_2, p_3) = (2 - \ffrac{1}{2 m_1} - \ts{\frac{1}{2
m_2}}) d^2 < 2 d^2.
\]

Assume now that $I_0 = \{ 1 \}$ and $I_1 = \{ 2, 3 \}$, thus $p_1
= 0$. Moreover, there exists $\lambda$ such that $p_2 =
\frac{m_2}{m_2 - 1} \lambda - d$ and $p_3 = \frac{m_3}{m_3 - 1}
\lambda - d$. Since $p + p_1 + p_2 + p_3 = d$, it follows that
$\lambda = \frac{3}{2 + \gamma}{d}$. The inequality $p_3 \leq
\frac{1}{m_3 - 1} d$ implies that $\lambda \leq d$, hence $\gamma
\geq 1$. By direct calculations we get
\[
f (p, p_1, p_2, p_3) = (2 - \ffrac{1}{2 m_1} - \ffrac{3 (\gamma -
1)}{2 (2 + \gamma)}) d^2 < 2 d^2.
\]

Let now $I_0 = \{ 1 \}$, $I_1 = \{ 2 \}$ and $I_3 = \{ 3 \}$. Then
$p_1 = 0$ and there exists $\lambda$ such that $p_2 =
\frac{m_2}{m_2 - 1} \lambda - d$ and $p_3 = d - \frac{m_3 - 2}{m_3
- 1} \lambda$. It follows that $\lambda = \frac{1}{\gamma} d$. The
inequality $p_3 > \frac{1}{m_3 - 1} d$ implies that $\lambda < d$,
hence $\gamma > 1$. One calculates that
\[
f (p, p_1, p_2, p_3) = (2 - \ffrac{1}{2 m_1} - \ffrac{\gamma -
1}{2 \gamma}) d^2 < d^2.
\]

We consider now the case $I_0 = \{ 1 \}$ and $I_2 = \{ 2, 3 \}$.
Then $p_1 = 0$ and there exists $\lambda$ such that $p_2 = d -
\frac{m_2 - 2}{m_2 - 1} \lambda$ and $p_3 = d - \frac{m_3 - 2}{m_3
- 1} \lambda$. It follows that $\lambda = \frac{1}{2 - \gamma} d$.
We get
\[
f (p, p_1, p_2, p_3) = (2  - \ffrac{\delta - 1}{2 (\delta - \gamma
+ 1)(2 - \gamma)}) d^2 \leq 2 d^2.
\]
Note that $\gamma < 2$ since in this case $m_2, m_3 > 2$. The
inequality is strict if and only if $\delta > 1$. Note that if
$\delta = 1$, then $\lambda = \frac{m_1 - 1}{m_1} d$ and $p_1 =
\frac{m_1}{m_1 - 1} \lambda - d$.

The next case is $I_1 = \{ 1, 2, 3 \}$. Then $p_1 = \frac{m_1}{m_1
- 1} \lambda - d$, $p_2 = \frac{m_2}{m_2 - 1} \lambda - d$ and
$p_3 = \frac{m_3}{m_3 - 1} \lambda - d$ for some $\lambda$. It
follows that $\lambda = \frac{4}{3 + \delta} d$. We get
\[
f (p,p_1, p_2, p_3) = (2 - \ffrac{2 (\delta - 1)}{3 + \delta}) d^2
\leq 2 d^2.
\]
Equality holds if and only if $\delta = 1$. If this is the case
then $\lambda = d$, $p_2 = d - \frac{m_2 - 2}{m_2 - 1} \lambda$
and $p_3 = d - \frac{m_3 - 2}{m_3 - 1} \lambda$.

Assume now that $I_1 = \{ 1, 2 \}$ and $I_2 = \{ 3 \}$. There
exists $\lambda$ such that $p_1 = \frac{m_1}{m_1 - 1} \lambda -
d$, $p_2 = \frac{m_2}{m_2 - 1} \lambda - d$ and $p_3 = d -
\frac{m_3 - 2}{m_3 - 1} \lambda$. It follows that $\lambda =
\frac{2}{1 + \delta} d$. The inequality $p_3 > \frac{1}{m_3 - 1}
d$ implies that $\lambda < d$, hence $\delta > 1$. We get
\[
f (p, p_1, p_2, p_3) = (2 - \ffrac{\delta - 1}{1 + \delta}) d^2 <
2 d^2.
\]

Finally, let $I_1 = \{ 1 \}$ and $I_2 = \{ 2, 3 \}$. There exists
$\lambda$ such that $p_1 = \frac{m_1}{m_1 - 1} \lambda - d$, $p_2
= d - \frac{m_2 - 2}{m_2 - 1} \lambda$ and $p_3 = d - \frac{m_3 -
2}{m_3 - 1} \lambda$. It follows that $(\delta - 1) \lambda = 0$.
If $\delta > 1$, then $\lambda = 0$ and $p_1 = - d < 0$, which is
impossible. Assume now that $\delta = 1$. The inequalities $p_1 >
0$ and $p_3 > \frac{1}{m_3 - 1} d$ imply that $\frac{m_1 - 1}{m_1}
d < \lambda < d$. One also checks that in this case \
\[
f (p, p_1, p_2, p_3) = 2 d^2,
\]
which finishes the proof.
\end{proof}

\subsection{} \label{subsectstrategy}
Our first aim is to prove Proposition~\ref{propineqnew} in the
following situation. Let $\bd \in \bR$ and $\bd' \in \bP$ be such
that $\bd' \neq 0$, $d_\alpha' = d_\alpha$, $d_\omega' = 0$, $\bd
- \bd' \in \bbN^{\Delta_0}$ and $p_{i, j}^{\bd} = 0$ for all $i
\in [1, n]$ and $j \in [1, m_i - 1]$. For simplicity we write in
this case $d$, $p$ and $p_1$, \ldots, $p_n$ instead of $d_\alpha$,
$p^{\bd}$ and $p_{1, m_1}^{\bd}$, \ldots, $p_{n, m_n}^{\bd}$,
respectively. Note that
\begin{multline*}
\langle \bd - \bd', \bd' \rangle = - d^2 - d p + \sum_{i \in [1,
n]} \Bigl(- d_{i, m_i - 1}' p_i + \sum_{j \in [1, m_i - 1]} (d_{i,
j - 1}' - d_{i, j}') d_{i, j}'\Bigr)
\end{multline*}
and $d = p + \sum_{i \in [1, n]} p_i$. Let $\delta_{i, j} = d_{i,
j - 1}' - d_{i, j}'$ for $i \in [1, n]$ and $j \in [1, m_i]$. Then
\[
\langle \bd - \bd', \bd' \rangle = - d^2 - d p + \sum_{i \in [1,
n]} S_i,
\]
where $S_i = - \delta_{i, m_i} p_i + \sum_{j < l \in [1, m_i]}
\delta_{i, j} \delta_{i, l}$. Note that $\delta_{i, j} \geq 0$ for
$i \in [1, n]$ and $j \in [1, m_i]$, $\sum_{j \in [1, n]}
\delta_{i, j} = d$ for $i \in [1, n]$, and $\delta_{i, 1} =
d_\alpha' - d_{i, 1}' \geq d_\alpha - d_{i, 1} = p_i$ for $i \in
[1, n]$.

Let $\frakO$ be the set of all pairs $(\bd, \bd')$ such that $\bd
\in \bR$, $p^{\bd} = 0$, $p_{i, j}^{\bd} = 0$ for $i \in [1, n]$,
$j \in [1, m_i - 1]$, $\bd' \in \bP$, $d_{i, j - 1}' > d_{i, j}'$
for $i \in [1, n]$ and $j \in [1, m_i - 1]$, $d_\omega' = 0$, and
$\bd - \bd' \in \bbN^{\Delta_0}$.

\subsection{} \label{specialcaseone}
Assume first that $n = 3$. It follows from the above paragraph and
Lemma~\ref{lemmineq} that
\[
\langle \bd - \bd', \bd' \rangle \leq - 2 d^2 + f (p, p_1, p_2,
p_3),
\]
where $f$ is as in~\ref{lemmfunf}. Lemma~\ref{lemmfunf} shows that
if $\frac{1}{m_1 - 1} + \frac{1}{m_2 - 1} + \frac{1}{m_3 - 1} \geq
1$, then
\[
\langle \bd - \bd', \bd' \rangle \leq 0,
\]
and if equality holds, then $\frac{1}{m_1 - 1} + \frac{1}{m_2 - 1}
+ \frac{1}{m_3 - 1} = 1$ and $p^{\bd} = 0$, which finishes the
proof of Proposition~\ref{propineqnew} in this case. Note also
that if equality holds then according to Lemmas~\ref{lemmineq}
and~\ref{lemmfunf}, $(\bd, \bd') \in \frakO$.

\subsection{} \label{specialcasetwo}
As the next case, consider $\bm = (2, 2, 2, m)$. It follows from
Lemma~\ref{lemmineq} that
\[
S_i \leq \ffrac{1}{4} d^2 - \ffrac{1}{2} d p_i + \ffrac{1}{4}
p_i^2
\]
for $i \in [1, 3]$, where we use notation introduced
in~\ref{subsectstrategy}. If $p_4 = d$, then $p = p_1 = p_2 = p_3
= 0$, and it follows from Corollary~\ref{coroineq} that $S_4 \leq
0$, hence
\[
\langle \bd - \bd', \bd' \rangle \leq - \ffrac{1}{4} d^2 < 0.
\]
On the other hand, if $d > p_4$, then using again
Corollary~\ref{coroineq},
\[
S_4 < \ffrac{1}{2} d^2 - \ffrac{1}{2} p_4^2,
\]
hence
\begin{multline*}
\langle \bd - \bd', \bd' \rangle
\\ %
< - \ffrac{3}{4} p^2 - \ffrac{1}{4} p_4^2 - p p_1 - p p_2 - p p_3
- \ffrac{1}{2} p p_4 - \ffrac{1}{2} p_1 p_2 - \ffrac{1}{2} p_1 p_3
- \ffrac{1}{2} p_2 p_3 \leq 0.
\end{multline*}

\subsection{} \label{specialcasethree}
Assume now that $\bm = (2, 2, 3, 3)$. It follows from
Lemma~\ref{lemmineq} that
\[
S_i \leq \ffrac{1}{4} d^2 - \ffrac{1}{2} d p_i + \ffrac{1}{4}
p_i^2
\]
for $i \in [1, 2]$. If $2 p_3 \leq d$ and $2 p_4 \leq d$, then
using again Lemma~\ref{lemmineq}, we get
\[
S_i \leq \ffrac{1}{3} d^2 - \ffrac{1}{3} d p_i + \ffrac{1}{3}
p_i^2
\]
for $i \in [3, 4]$, hence
\begin{multline*}
\langle \bd - \bd', \bd' \rangle \leq - \ffrac{5}{6} p^2 -
\ffrac{1}{12} p_1^2 - \ffrac{1}{12} p_2^2 - \ffrac{7}{6} p p_1 -
\ffrac{7}{6} p p_2 - \ffrac{5}{6} p p_3 - \ffrac{5}{6} p p_4
\\ %
- \ffrac{2}{3} p_1 p_2 - \ffrac{1}{3} p_1 p_3 - \ffrac{1}{3} p_1
p_4 - \ffrac{1}{3} p_2 p_3 - \ffrac{1}{3} p_2 p_4 - \ffrac{1}{6}
p_3 (d - 2 p_3) - \ffrac{1}{6} p_4 (d - 2 p_4) \leq 0.
\end{multline*}
Moreover, if equality holds, then $p = p_1 = p_2 = 0$ and $p_3 =
p_4$. Applying in addition once more Lemma~\ref{lemmineq} we get
that if equality holds then $(\bd, \bd') \in \frakO$.

As the next case, consider $2 p_4 > d$, i.e., $p_4 > p + p_1 + p_2
+ p_3$. Then in particular $2 p_3 \leq d$, hence
\[
S_3 \leq \ffrac{1}{3} d^2 - \ffrac{1}{3} d p_3 + \ffrac{1}{3}
p_3^2, \qquad S_4 \leq d p_4 - p_4^2,
\]
and
\begin{multline*}
\langle \bd - \bd', \bd' \rangle \leq -\ffrac{7}{6} p^2 -
\ffrac{5}{12} p_1^2 - \ffrac{5}{12} p_2^2 - \ffrac{1}{6} p_3^2 -
\ffrac{1}{6} p_4^2  - \ffrac{11}{6} p p_1 - \ffrac{11}{6} p p_2
\\ %
- \ffrac{5}{3} p p_3 - \ffrac{1}{3} p p_4 - \ffrac{4}{3} p_1 p_2 -
\ffrac{7}{6} p_1 p_3 + \ffrac{1}{6} p_1 p_4 - \ffrac{7}{6} p_2 p_3
+ \ffrac{1}{6} p_2 p_4 + \ffrac{1}{3} p_3 p_4.
\end{multline*}
One easily checks that the above expression is decreasing when
considered as a function of $p_4$ for $p_4 > p + p_1 + p_2 + p_3$.
Moreover, for $p_4 = p + p_1 + p_2 + p_3$ we get
\[
- \ffrac{5}{3} p^2 - \ffrac{5}{12} p_1^2 - \ffrac{5}{12} p_2^2 -
\ffrac{7}{3} p p_1 - \ffrac{7}{3} p p_2 - 2 p p_3 - \ffrac{4}{3}
p_1 p_2 - p_1 p_3 - p_2 p_3,
\]
hence $\langle \bd - \bd', \bd \rangle < 0$ in this case.

\subsection{} \label{specialcasefour}
The final case we have to consider is $\bm = (2, 2, 2, 2, 2)$. Let
$\frakO'$ be the set of all pairs $(\bd, \bd') \in \frakO$ such
that $\bd = d \be_{i, 2}$ for a positive integer $d$ and some $i
\in [1, 5]$, $d_\alpha' \in [1, d]$, $d_{i, 1}' = 0$, $d_{j, 1}' =
\frac{1}{2} d_\alpha'$ for $j \in [1, 5]$, $j \neq i$, and
$d_\omega' = 0$.

Using Lemma~\ref{lemmineq} we get that
\[
S_i \leq \ffrac{1}{4} d^2 - \ffrac{1}{2} d p_i + \ffrac{1}{4}
p_i^2
\]
for $i \in [1, 5]$, hence
\[
\langle \bd - \bd', \bd' \rangle = - \ffrac{3}{4} p^2 - \sum_{i
\in [1, 5]} p p_i - \ffrac{1}{2} \sum_{i < j \in [1, 5]} p_i p_j
\leq 0.
\]
Moreover, if equality holds, then $p = 0$ and there exists $i \in
[1, 5]$ such that $p_i = d$ and $p_j = 0$ for $j \in [1, 5]$, $j
\neq i$. Finally, it follows from Lemma~\ref{lemmineq} that in the
case of equality $(\bd, \bd') \in \frakO'$.

\subsection{} \label{lemmredone}
We show now that we can reduce the proof of
Proposition~\ref{propineqnew} to the special situation considered
in the previous paragraphs. We first show that we may assume that
$d_\omega' = 0$.

\begin{lemm}
Let $\bd \in \bR$ and $\bd' \in \bP$ be such that $\bd' \neq 0$
and  $\bd - \bd' \in \bR + \bQ$. If $d_\omega' > 0$, then $\bd' -
\bh \in \bP$, $\bd' - \bh \neq 0$,  $\bd - (\bd' - \bh) \in \bR +
\bQ$, and
\[
\langle \bd - \bd', \bd' \rangle = \langle \bd - (\bd' - \bh),
\bd'- \bh \rangle.
\]
\end{lemm}

\begin{proof}
The former three assertions are obvious ($\bd' - \bh \neq 0$,
since $\bh \not \in \bP$), the latter follows by direct
calculations.
\end{proof}

\subsection{} \label{lemmredtwo}
The second reduction is the following.

\begin{lemm}
Let $\bd \in \bR$ and $\bd' \in \bP$ be such that $\bd' \neq 0$,
$\bd - \bd' \in \bR + \bQ$ and $d_\omega' = 0$. If $i \in [1, n]$
and $j \in [1, m_i - 1]$ are such that $p_{i, j}^{\bd} > 0$ and
$p_{i, j}^{\bd - \bd'} > 0$, then $\bd - \be_{i, j} \in \bR$,
$(\bd - \be_{i, j}) - \bd' \in \bR + \bQ$ and
\[
\langle \bd - \bd', \bd' \rangle \leq \langle (\bd - \be_{i, j}) -
\bd', \bd' \rangle.
\]
Moreover, if $(\bd - \be_{i, j}, \bd') \in \frakO$, then the above
inequality is strict.
\end{lemm}

\begin{proof}
Obviously, $p_{i, j}^{\bd} > 0$ implies that $\bd - \be_{i, j} \in
\bR$. Similarly, $p_{i, j}^{\bd - \bd'} > 0$ implies that $(\bd -
\be_{i, j}) - \bd' = (\bd - \bd') - \be_{i, j} \in \bR + \bQ$.
Moreover,
\[
\langle (\bd - \be_{i, j}) - \bd', \bd' \rangle - \langle \bd -
\bd', \bd' \rangle = d_{i, j - 1}' - d_{i, j}' \geq 0.
\]
Finally, if $(\bd - \be_{i, j}, \bd') \in \frakO$, then $d_{i, j -
1}' > d_{i, j}'$, hence the above inequality is strict.
\end{proof}

\subsection{} \label{lemmredtwobis}
A more complicated version of the above reduction is the
following.

\begin{lemm}
Let $\bd \in \bR$ and $\bd' \in \bP$ be such that $\bd' \neq 0$,
$\bd - \bd' \in \bR + \bQ$ and $d_\omega' = 0$. If $i \in [1, n]$
and $j \in [1, m_i - 1]$ are such that $p_{i, j}^{\bd} > 0$,
$p_{i, l}^{\bd} = 0$ for all $l \in [j + 1, m_i - 1]$, and $p_{i,
j}^{\bd - \bd'} = 0$, then $\bd - \be_{i, j} \in \bR$, $\bd' -
\be_{i, j} \in \bP$, $\bd' - \be_{i, j} \neq 0$, $(\bd - \be_{i,
j}) - (\bd' - \be_{i, j}) \in \bR + \bQ$ and
\[
\langle \bd - \bd', \bd' \rangle \leq \langle (\bd - \be_{i, j}) -
(\bd' - \be_{i, j}), \bd' - \be_{i, j} \rangle.
\]
Moreover, if $(\bd - \be_{i, j}, \bd' - \be_{i, j}) \in \frakO$,
then the above inequality is strict.
\end{lemm}

\begin{proof}
Obviously, $(\bd - \be_{i, j}) - (\bd' - \be_{i, j}) = \bd - \bd'
\in \bR + \bQ$ and $p_{i, j}^{\bd} > 0$ implies that $\bd -
\be_{i, j} \in \bR$. Moreover, $\bd' - \be_{i, j} \neq 0$, since
$\be_{i, j} \not \in \bP$. Note that
\begin{align*}
d_{i, l} - d_{i, l}' & = p^{\bd - \bd'} + \sum_{\substack{s \in
[1, n] \\ s \neq i}} p_{s, m_s}^{\bd - \bd'} + p_{i, l}^{\bd -
\bd'}
\\ %
\intertext{for $l \in [1, m_i - 1]$, and} %
d_{i, m_i} - d_{i, m_i}' & = p^{\bd - \bd'} + \sum_{s \in [1, n]}
p_{s, m_s}^{\bd - \bd'} + p_\omega^{\bd - \bd'}.
\end{align*}
Thus our assumption implies that
\begin{equation} \label{eqmin}
d_{i, j} - d_{i, j}' = \min_{l \in [1, m_i]} (d_{i, l} - d_{i,
l}').
\end{equation}
In particular,
\begin{equation} \label{ineqjj}
d_{i, j} - d_{i, j}' \leq d_{i, j + 1} - d_{i, j + 1}',
\end{equation}
and the above inequality is strict if $j = m_i - 1$ (note that
$\bd' \neq 0$ implies that $p_\omega^{\bd - \bd'} \neq 0$).

We show now that $\bd' - \be_{i, j} \in \bP$. In order to do this
we have to prove that $d_{i, j}' - d_{i, j + 1}' > 0$. If $j < m_i
- 1$, then using that $p_{i, j + 1}^{\bd} = 0$ and~\eqref{ineqjj}
we get
\[
d_{i, j}' - d_{i, j + 1}' \geq d_{i, j} - d_{i, j + 1} = p_{i,
j}^{\bd} - p_{i, j + 1}^{\bd} = p_{i, j}^{\bd} > 0.
\]
If $j = m_i - 1$, then $d_{i, j + 1}' = d_\omega' = 0$, so we have
to prove that $d_{i, j}' > 0$. Choose $l \in [1, m_i]$ such that
$p_{i, l}^{\bd} = 0$. It follows similarly as above that $d_{i, l}
< d_{i, j}$. Using~\eqref{eqmin} we get
\[
d_{i, j}' > d_{i, l}' \geq 0,
\]
thus the claim follows.

In order to prove the required inequality note that
\begin{multline*}
\langle (\bd - \be_{i, j}) - (\bd' - \be_{i, j}), \bd' - \be_{i,
j} \rangle - \langle \bd - \bd', \bd' \rangle
\\ %
= (d_{i, j + 1} - d_{i, j + 1}') - (d_{i, j} - d_{i, j}'),
\end{multline*}
hence the claim follows from~\eqref{ineqjj}. It also follows that
if $j = m_i - 1$, then the inequality is strict. Finally assume
that $j \in [1, m_i - 2]$ and $(\bd - \be_{i, j}, \bd' - \be_{i,
j}') \in \frakO$. This implies that $d_{i, j} = d_{i, j + 1} + 1$
and $d_{i, j}' \geq d_{i, j + 1}' + 2$, which finishes the proof.
\end{proof}

\subsection{} \label{lemmredthree}
The last reduction is the following.

\begin{lemm}
Let $\bd \in \bR$ and $\bd' \in \bP$ be such that $\bd' \neq 0$,
$\bd - \bd' \in \bR + \bQ$, $d_\omega' = 0$, and $p_{i, j}^{\bd} =
0$ for all $i \in [1, n]$ and $j \in [1, m_i - 1]$. If $d_\alpha'
< d_\alpha$, then there exists $i \in [1, n]$ such that $\bd -
\be_{i, m_i} \in \bR$, $(\bd - \be_{i, m_i}) - \bd' \in \bR + \bQ$
and
\begin{equation} \label{ineq3}
\langle \bd - \bd', \bd' \rangle \leq \langle (\bd - \be_{i, m_i})
- \bd', \bd' \rangle.
\end{equation}
Moreover, if $(\bd - \be_{i, m_i}, \bd') \in \frakO$, then $(\bd,
\bd') \in \frakO$. Finally, if $\bm = (2, 2, 2, 2, 2)$, $(\bd -
\be_{i, m_i}, \bd') \in \frakO'$ and $(\bd, \bd') \not \in
\frakO'$, then the above inequality is strict.
\end{lemm}

\begin{proof}
We first show the existence of $i \in [1, n]$ such that $\bd -
\be_{i, m_i} \in \bR$ and $(\bd - \be_{i, m_i}) - \bd' \in \bR +
\bQ$. Observe that $0 < d_\alpha' < d_\alpha = p^{\bd} + \sum_{i
\in [1, n]} p_{i, m_i}^{\bd}$, hence either $p^{\bd} > 0$ or there
exists $i \in [1, n]$ such that $p_{i, m_i}^{\bd} > 0$. Similarly,
since $d_\alpha - d_\alpha' > 0$, either $p^{\bd - \bd'} > 0$ or
there exists $i \in [1, n]$ such that $p_{i, m_i}^{\bd - \bd'} >
0$. Note that if $p^{\bd} > 0$, then $\bd - \be_{i, m_i} \in \bR$
for all $i \in [1, n]$, since $\bh - \be_{i, m_i} = \sum_{j \in
[1, m_i - 1]} \be_{i, j}$. Again similarly, if $p^{\bd - \bd'} >
0$, then $(\bd - \be_{i, m_i}) - \bd' = (\bd - \bd') - \be_{i,
m_i} \in \bR + \bQ$ for all $i \in [1, n]$. Thus it remains to
show that, if $p^{\bd} = 0 = p^{\bd - \bd'}$, then there exists $i
\in [1, n]$ such that $p_{i, m_i}^{\bd}, p_{i, m_i}^{\bd - \bd'} >
0$. Without loss of generality we may assume that $p_{1,
m_1}^{\bd}, \ldots, p_{s, m_s}^{\bd} > 0$ and $p_{s + 1, m_{s +
1}}^{\bd} = \cdots = p_{n, m_n}^{\bd} = 0$ for some $s \in [1,
n]$. Then for $i \in [s + 1, n]$ and $j \in [1, m_i - 1]$, $d_{i,
j} = d_\alpha$ and $(d_{i, j} - d_{i, j}') - (d_\alpha -
d_\alpha') = p_{i, j}^{\bd - \bd'} - p_{i, m_i}^{\bd - \bd'}$,
hence $p_{i, m_i}^{\bd - \bd'} - p_{i, j}^{\bd - \bd'} = d_{i, j}'
- d_\alpha' \leq 0$. Consequently, $p_{i, m_i}^{\bd - \bd'} = \min
\{ p_{i, j}^{\bd - \bd'} \mid j \in [1, m_i] \} = 0$ for $i \in [s
+ 1, n]$. Since $d_\alpha - d_\alpha' > 0$, it follows that there
exists $i \in [1, s]$ such that $p_{i, m_i}^{\bd - \bd'} > 0$.

Note that
\[
\langle (\bd - \be_{i, m_i}) - \bd', \bd' \rangle - \langle \bd -
\bd', \bd' \rangle = d_{i, m_i - 1}' \geq 0.
\]
Obviously, if $(\bd - \be_{i, m_i}, \bd') \in \frakO$, then $(\bd,
\bd') \in \frakO$. Finally, assume that $\bm = (2, 2, 2, 2, 2)$
and $(\bd - \be_{i, m_i}, \bd') \in \frakO'$. In particular, $\bd
- \be_{i, m_i} = d \be_{j, m_j}$ for a positive integer $d$ and $j
\in [1, n]$. If $(\bd, \bd') \not \in \frakO'$, then $j \neq i$,
hence $d_{i, m_i - 1}' = \ffrac{1}{2} d \neq 0$, and the above
inequality is strict, which finishes the proof.
\end{proof}

\subsection{}
We can complete now the proof of Proposition~\ref{propineqnew}.
Let $\Lambda$ be a canonical algebra of type $\bm$, $\bd \in \bR$
and $\bd' \in \bP$ be such that $\bd' \neq 0$ and $\bd - \bd' \in
\bR + \bQ$. It follows from Lemma~\ref{lemmredone}, that we may
assume $d_\omega' = 0$. It follows by an easy induction that there
exists a sequence $(\bd^{(s)}, \bd'^{(s)})$, $s \in [0, l]$, such
that $\bd^{(0)} = \bd$, $\bd'^{(0)} = \bd'$, $(\bd^{(s)},
\bd'^{(s)})$ is obtained from $(\bd^{(s - 1)}, \bd'^{(l - s)})$,
$s \in [1, l]$, by applying one of the reductions described in
Lemmas~\ref{lemmredtwo}--\ref{lemmredthree}, $d_\alpha^{(l)} =
d_\alpha'^{(l)}$, and $p_{i, j}^{\bd^{(l)}} = 0$ for all $i \in
[1, n]$ and $j \in [1, m_i - 1]$. In particular we know that
\[
\langle \bd - \bd', \bd' \rangle \leq \langle \bd^{(l)} -
\bd'^{(l)}, \bd'^{(l)} \rangle \leq 0,
\]
where the latter inequality follows
from~\ref{specialcaseone}--\ref{specialcasefour}. Moreover, the
latter inequality is strict if $\sum_{i \in [1, n]} \frac{1}{m_i -
1} > 2 n - 5$

Now assume that $\sum_{i \in [1, n]} \frac{1}{m_i - 1} = 2 n - 5$
and $\langle \bd - \bd', \bd' \rangle = 0$. Then $\langle
\bd^{(l)} - \bd'^{(l)}, \bd'^{(l)} \rangle = 0$ and consequently
\[
\langle \bd^{(s - 1)} - \bd'^{(s - 1)}, \bd'^{(s - 1)} \rangle =
\langle \bd^{(s)} - \bd'^{(s)}, \bd'^{(s)} \rangle
\]
for all $s \in [1, l]$. It also follows from \ref{specialcaseone},
\ref{specialcasethree} and \ref{specialcasefour} that $(\bd^{(l)},
\bd'^{(l)}) \in \frakO$, hence using
Lemmas~\ref{lemmredtwo}--\ref{lemmredthree} we get by induction
that for all $s \in [0, l]$, $(\bd^{(s)}, \bd'^{(s)}) \in \frakO$.
In particular, $p^{\bd} = p^{\bd^{(0)}} = 0$. With similar
arguments we prove that $(\bd, \bd') \in \frakO'$ if $\bm = (2, 2,
2, 2, 2)$ and $\langle \bd - \bd', \bd' \rangle = 0$, which
implies the last assertion of Proposition~\ref{propineqnew}.

\makeatletter
\def\@secnumfont{\mdseries} 
\makeatother

\section{Counterexamples} \label{sectexamp}

In this section we present for a canonical algebra of type $\bm$
such that
\[
\sum_{i \in [1, n]} \ffrac{1}{m_i - 1} = 2 n - 5 \text{ (} < 2 n -
5 \text{, respectively)},
\]
examples of dimension vectors $\bd' \in \bP$ and $\bd'' \in \bQ$
such that $\bd' + \bd'' \in \bR$ and
\[
\langle \bd'', \bd' \rangle = 0 \text{ (} > 0 \text{,
respectively)}.
\]
Together with Propositions~\ref{criterioncomplete},
\ref{criterion}, \ref{propnorm} and~\ref{propineqnew}, it will
finish the proof of Theorems~\ref{theo1}, \ref{theo2}
and~\ref{theo3}.

\subsection{}
Let $\Lambda$ be a canonical algebra of type $(m_1, m_2, m_3)$
such that
\[
\delta = \ffrac{1}{m_1 - 1} + \ffrac{1}{m_2 - 1} + \ffrac{1}{m_3 -
1} \leq 1.
\]
Note that our assumption implies that $m_1, m_2, m_3 > 2$. Let
\[
m = (m_1 - 1) (m_2 - 1) (m_3 - 1) (m_1 - 2) (m_2 - 2) (m_3 - 2).
\]
Define $\bd'$ and $\bd''$ by
\begin{align*}
d_\alpha' & = \delta m, & d_\alpha'' & = 0,
\\ %
d_{i, j}' & = \ffrac{(\delta (m_i - 1) - 1) (m_i - j - 1)}{(m_i -
1) (m_i - 2)} m, & d_{i, j}'' & = \ffrac{(\delta (m_i - 1) - 1 )
(j - 1)}{(m_i - 1) (m_i - 2)} m,
\\ %
& & & \qquad i \in [1, 3], \, j \in [1, m_i - 1],
\\ %
d_\omega' & = 0, & d_\omega'' & = \delta m.
\end{align*}
Then $\bd' \in \bP$, $\bd'' \in \bQ$,
\begin{align*}
\bd' + \bd'' & = \ffrac{m}{m_1 - 1} \be_{1, m_1} + \ffrac{m}{m_2 -
1} \be_{2, m_2} + \ffrac{m}{m_3 - 1} \be_{3, m_3} \in \bR
\\ %
\intertext{and} %
\langle \bd'', \bd' \rangle & = \sum_{i \in [1, 3]} \sum_{j \in
[2, m_i - 1]} d_{i, j}'' (d_{i, j}' - d_{i, j - 1}') + d_\omega''
d_\alpha'
\\ %
& = \ffrac{1}{2} \bigl( - \sum_{i \in [1, 3]} \ffrac{(\delta (m_i
- 1) - 1)^2}{(m_i - 1) (m_i - 2)} + 2 \delta^2 \bigr) m^2
\\ %
& = \ffrac{1}{2} (- \delta^2 - \delta^2 \delta' + 2 \delta \delta'
+ \delta - \delta') m^2 = \ffrac{1}{2} (1 - \delta) (\delta
\delta' + \delta - \delta')
\\ %
& = \ffrac{1}{2} (1 - \delta) \bigl( \sum_{i \neq j \in [1, 3]}
\ffrac{1}{(m_i - 1) (m_j - 2)} \bigr) m^2 \geq 0,
\intertext{where} %
\delta' & = \ffrac{1}{m_1 - 2} + \ffrac{1}{m_2 - 2} +
\ffrac{1}{m_3 - 2}.
\end{align*}
The above inequality is strict if $\delta < 1$.

\subsection{}
Let $\Lambda$ be a canonical algebra of type $(m_1, m_2, m_3,
m_4)$ such that
\[
\ffrac{1}{m_1 - 1} + \ffrac{1}{m_2 - 1} + \ffrac{1}{m_3 - 1} +
\ffrac{1}{m_4 - 1} \leq 3.
\]
The above assumption implies in particular that, without loss of
generality, we may assume that $m_3, m_4 > 2$. Let
\[
m = 2 m_1 m_2 (m_3 - 2) (m_4 - 2).
\]
Define $\bd'$ and $\bd''$ by
\begin{align*}
d_\alpha' & = m, & d_\alpha'' & = 0,
\\ %
d_{i, j}' & = \ffrac{m_i - j}{m_i} m, & d_{i, j}'' & =
\ffrac{j}{m_i} m, & & i \in [1, 2], \, j \in [1, m_i - 1],
\\ %
d_{i, j}' & = \ffrac{m_i - j - 1}{2 (m_i - 2)} m, & d_{i, j}'' & =
\ffrac{j - 1}{2 (m_i - 2)} m, & & i \in [3, 4], \, j \in [1, m_i -
1],
\\ %
d_\omega' & = 0, & d_\omega'' & = m.
\end{align*}
Then $\bd' \in \bP$, $\bd'' \in \bQ$,
\begin{align*}
\bd' + \bd'' & = \ffrac{m}{2} \be_{3, m_3} + \ffrac{m}{2} \be_{3,
m_4} \in \bR
\\ %
\intertext{and} %
\langle \bd'', \bd' \rangle & = (\ffrac{3}{4} - \ffrac{1}{2 m_1} -
\ffrac{1}{2 m_2} - \ffrac{1}{8 (m_3 - 2)} - \ffrac{1}{8 (m_4 -
2)}) m^2 \geq 0.
\end{align*}
The inequality is strict if $\bm \neq (2, 2, 3, 3)$.

\subsection{}
Let $\Lambda$ be a canonical algebra of type $(m_1, \ldots, m_n)$
for $n \geq 5$. We may assume, without loss of generality, that
$m_n = \min (m_1, \ldots, m_n)$. Let
\[
m = m_1 \cdots m_{n - 1}.
\]
Define $\bd'$ and $\bd''$ by
\begin{align*}
d_\alpha' & = m, & d_\alpha'' & = 0,
\\ %
d_{i, j}' & = \ffrac{m_i - j}{m_i} m, & d_{i, j}'' & =
\ffrac{j}{m_i} m, & & i \in [1, n - 1], \, j \in [1, m_i - 1],
\\ %
d_{n, j}' & = 0, & d_{n, j}'' & = 0, & & j \in [1, m_n - 1],
\\ %
d_\omega' & = 0, & d_\omega'' & = m.
\end{align*}
Then $\bd' \in \bP$, $\bd'' \in \bQ$,
\begin{align*}
\bd' + \bd'' & = m \be_{n, m_n}
\\ %
\intertext{and} %
\langle \bd'', \bd' \rangle & = \ffrac{1}{2} \bigl( n - 3 -
\sum_{i \in [1, n - 1]} \ffrac{1}{m_i} \bigr) m^2 \geq 0.
\end{align*}
The inequality is strict if $\bm \neq (2, 2, 2, 2, 2)$ (remember,
that $m_n = \min (m_1, \ldots, m_n)$).

\subsection{}
Note that in all the above examples $p^{\bd} = 0$ for $\bd = \bd'
+ \bd''$. Moreover, $\bd$ is not sincere for $n \geq 5$. In order
to complete the proof of Theorem~\ref{theo3} we have to present
examples with $p^{\bd} > 0$ and $\langle \bd'', \bd' \rangle > 0$,
for canonical algebras $\Lambda$ of type $(m_1, \ldots, m_n)$ with
$\sum_{i \in [1, n]} \frac{1}{m_i - 1} < 2 n - 5$. It will also
complete the proof of Theorem~\ref{theo2}, since $\bd \in \bR$
with $p^{\bd} > 0$ is sincere.

Let $\Lambda$ be an algebra of the above form. It follows from the
preceding paragraphs that there exist dimension vectors $\bd' \in
\bP$ and $\bd'' \in \bQ$ such that $\bd' + \bd'' \in \bR$ and
$\langle \bd'', \bd' \rangle > 0$. Choose a positive integer $q$
such that
\[
q \langle \bd'', \bd' \rangle + \langle \bd'', \bh \rangle
> 0.
\]
Then $\hat{\bd}' = q \bd' + \bh \in \bP$, $\hat{\bd}'' = q \bd''
\in \bQ$, $\bd = \hat{\bd}' + \hat{\bd}'' = \bh + q (\bd' + \bd'')
\in \bR$, $p^{\bd} > 0$ and
\[
\langle \hat{\bd}'', \hat{\bd}' \rangle = q^2 \langle \bd'', \bd'
\rangle + q \langle \bd'', \bh \rangle
> 0.
\]

\bibsection

\end{document}